\documentclass[11pt]{article}
\newfont{\inhead}{eufm10}

\renewcommand{\u}{{\mbox{\inhead{u}}}}

\newcommand{\pbar}{{\bar\partial}}

\newcommand{\Z}{\mathbb{Z}}
\newcommand{\R}{\mathbb{R}}
\newcommand{\C}{\mathbb{C}}

\newcommand{\Q}{\mathbb{Q}}

\newcommand{\Aa}{{\cal A}}   

\newcommand{\Dd}{{\cal D}}

\newcommand{\Ff}{{\cal F}}

\newcommand{\Jj}{{\cal J}}
\newcommand{\Mm}{{\cal M}}   

\newcommand{\Oo}{{\cal O}}

\newcommand{\Ss}{{\cal S}}

\newcommand{\Tt}{{\cal T}}

\newcommand{\coker}{{\rm coker }}  
\newcommand{\im}{{\rm im }}        

\newcommand{\Diff}{{\rm Diff}}        
\newcommand{\Symp}{{\rm Symp}}        
\newcommand{\Hom}{{\rm Hom}}          

\newcommand{\Teich}{{\rm Teich}}
\newcommand{\Smooth}{{\rm Smooth}}

\def\NABLA#1{{\mathop{\nabla\kern-.5ex\lower1ex\hbox{$#1$}}}}
\def\Nabla#1{\nabla\kern-.5ex{}_#1}
\newcommand{\Tilde}{\widetilde}

\newcommand{\QED}{\hfill$\square$\medskip}

\newtheorem{theorem}{Theorem}[section]
\newtheorem{corollary}[theorem]{Corollary}

\newtheorem{lemma}[theorem]{Lemma}
\newtheorem{proposition}[theorem]{Proposition}
\newtheorem{definition}[theorem]{Definition}
\newtheorem{remark}[theorem]{Remark}

\usepackage{amssymb}
\usepackage{amsmath}
\usepackage[all]{xy}

\begin{document}

\title{Whitehead products in symplectomorphism groups and Gromov-Witten invariants} 
\author{Olgu\c ta Bu\c se}
\date{\today}


\maketitle

\abstract{Consider any symplectic ruled surface $(M^g_{\lambda},\omega_{\lambda})=
          (\Sigma_g \times S^2, \lambda \sigma_{\Sigma_g} \oplus \sigma_{S^2})$. We
          compute all natural equivariant Gromov-Witten invariants
          $EGW_{g,0}(M^g_{\lambda};H_k, A-kF)$ for all hamiltonian circle
          actions $H_k$ on $M^g_{\lambda}$, where $A=[\Sigma_g \times pt]$ and
          $F= [pt \times S^2]$. We use these invariants to show the
          nontriviality of certain higher order Whitehead products that live
          in the homotopy groups of the symplectomorphism groups
          $G_{\lambda}^g$, $g \geq 0$. Our results are sharper when $g=0,1$ and enable us
          to answer a question posed by D.McDuff in \cite{M} in the case $g=1$
          and provide a new interpretation of the multiplicative structure in
          the ring $H^*(BG^0_{\lambda} ;\Q)$ studied by Abreu-McDuff in
          \cite{AM}.}

\section{Introduction}

We study some top
ological aspects of symplectomorphism groups $G_{\lambda}$ of a continuous family of symplectic structures  $(M,\omega_{\lambda}) $, $\lambda \geq 0$ on a given compact manifold. 

In section two we provide preliminary definitions
and results on symplectic fibrations, topological aspects of symplectomorphism
groups and Whitehead (and hence Samelson) products. Higher order Whitehead products of elements in  $\pi_*BG_{\lambda}$ are subsets of some $\pi_NBG_{\lambda}$ and measure the obstructions of extending existing maps defined on the codimension $2$ skeleton of an appropriate product of spheres to the product itself. These sets are nonempty only when the lower order products contain the nullhomotopic class. Samelson products  in $G_{\lambda}$ are desuspensions of Whitehead products.

We make use of these obstruction theoretic properties in  section three where we introduce and prove proposition-construction \ref{Whitehead}. This proposition relates, in certain circumstances, higher order Whitehead products to towers of symplectic fibrations over projective spaces.

Section four provides background material on parametric Gromov-Witten invariants PGW and provides a relation between these invariants and Whitehead products in the
symplectomorphism groups. Roughly speaking, one can relate continuous deformations with respect to the parameter $\lambda$ of Whitehead products in $BG_{\lambda}$ to a symplectic deformation problem on symplectic fibrations. Parametric Gromov -Witten invariants are precisely invariants of such fibrations.

Computing PGW invariants is, aside from their consequences on the symplectomorphism groups, of interest on its own. 
We provide such computation in section five where we study  ruled surfaces $(M^g_{\lambda},\omega_{\lambda})= (\Sigma_g \times S^2, \lambda \sigma_{\Sigma_g} \oplus \sigma_{S^2})$, with symplectomorphism groups $G^g_{\lambda}$.
If $k =\lfloor \lambda \rfloor$  then $M^g_{\lambda}$  admits $k$ different hamiltonian circle actions $H_i$ each with two fixed point sets given by holomorphic curves in classes $A \pm i F$,$1 \leq i \leq k$. Equivariant Gromov-Witten invariants EGW count precisely $H_i$ invariant curves. We show that the only {\it natural} EGW (those counting generically isolated curves with no marked points in certain associated fibrations) are given by:

\begin{theorem}\label{egwruled}
  For any arbitrary genus $g$, and a hamiltonian circle action with Lie group
  $H_k$ on $M_{\lambda}^g, \lambda > k$ as in (\ref{actionruledsymp}), and
  equivariant almost complex structure $J^{(k),g}$ we have
  $$
    EGW^{J^{(k),g} } _{g,0}(M^g _{\lambda};H_k; s_{A-kF})=
    \pm 1 \cdot \u^{2k+g-1} \in H^*( B S^1, \Q)
  $$
\end{theorem}

In order to successfully tie our result on equivariant Gromov-Witten invariants to nontrivial Samelson products in symplectomorphism groups of ruled surfaces we need to make use of several previously known results. Building on results from M. Gromov \cite{G}, M. Abreu \cite{A} and M. Abreu- D. McDuff \cite{AM} these results are provided by D. McDuff in  \cite{M}. Denote by  $\Aa_{\lambda}^g$ the space of almost complex structures $J$ that tame some symplectic form cohomologous to $\omega_{\lambda}$.
The results in \cite{M}
relate the homotopy groups of $G_{\lambda}^g$ to a good understanding of strata  $\Aa_{\lambda,k }^g$ of almost complex structures that admit curves in class $A-kF$. A complete description of this stratification, for $g>0$, is prevented by difficult gluing problems. When seeking nontrivial Samelson products  in $\pi_*G_{\lambda}^g$, we circumvent some of these issues by using the topological structures laid out in the previous sections, and obtain:

\begin{proposition}\label{fragileruledhg}
  
  {\renewcommand{\labelenumi}{{\bf (\roman{enumi})}}
  \begin{enumerate}
    \item There exists an element $\widetilde{\gamma} \in \pi_1 G_1^1 $ such
          that $[\widetilde{\gamma},\widetilde{\gamma}]_s \in \pi_2 G_{1}^1 $
          is a nontrivial element that disappears in $\pi_2G_{\lambda}^1$,
          $\lambda>1$.
    \item For all $k \geq 1 $ and $ k <\lambda \leq k+1 $ there exist elements
          $ \widetilde{\gamma}^0_{\lambda} \in \pi_1 G_{\lambda}^0 \otimes \Q$ such that
          the Samelson product of order $2 k +1$, $S^{(2k+1)}
          (\widetilde{\gamma}^0_{\lambda}) = \{ 0, \widetilde{w_k} \}
          \subset \pi_{4k} (G^0_{k+1}) $ where $\widetilde{w_k}$ is a
          nontrivial homotopy class that disappears when $\lambda >k+1$.
    \item For all genus $g \geq 2$ and all $k > [g/2 ] $ there exist elements
          $\widetilde{\gamma}_{\lambda}^g\in \pi_1G^g_k \otimes \Q$ and nonvanishing
          Samelson product of order $p$ with $g \leq p \leq 2 k +g -1$, $0
          \neq  \widetilde{w}_p^g \in S^{(p)} (\widetilde{\gamma}^g_{\lambda}) \subset \pi_{2p-2}
          (G_{k}^g)$ .
  \end{enumerate} }
\end{proposition}

In the case when $g=0$ we use this theorem and the additive structure on
$\pi_*G^0_{\lambda} \otimes \Q$ from \cite{AM}, to give a new proof for the
following:

\begin{theorem}{(Abreu-McDuff)\cite{AM}}\label{rationaltypeclass}
  Fix an integer $k\geq0$. For $k <\lambda \leq k+1$ we have
  \begin{equation}\label{classring}
    H^*(BG_{\lambda}^0, \Q) = S (A,X,Y) / \{ A(X-Y)(X -4Y) \cdots ( X-k^2Y)
    =0 \}
  \end{equation}
  with $deg A=2$ and $deg X =deg Y= 4$
\end{theorem}

\section{Preliminaries}

\subsection{Symplectic fibrations}\label{sect:symplecticfibration} 
Consider a triple $(M,\omega_0,J)$ where $J$ is an almost complex structure
that tames $\omega_0$ and has a canonical class $c_1(M) $.

\begin{definition}\label{def:Fibration}
  A locally trivial fibration $ \pi: Q \rightarrow B$ is a {\bf symplectic
  fibration} if the fiber is the compact symplectic manifold $(M,\omega_0)$
  and there exist a two form $\Lambda_0$ on $Q$ which is vertically closed
  i.e. $i(v_1,v_2)d\Lambda_0=0$ for all vertical vectors $v_i$ and whose
  restriction to each fiber is the symplectic form of the fiber.
\end{definition}

As shown in \cite{GS}, such forms correspond to symplectic connections on the
fibration.  Consider $(U_{\alpha})$ an atlas covering the base $B$ and a
trivialization $\phi_{\alpha}: \pi^{-1} (U_{\alpha}) \rightarrow M \times
U_{\alpha} $, that yields a collection of transition maps $\phi_{\alpha \beta}
:U_{\alpha} \cup U_{ \beta} \rightarrow\rm{ Diff} (M)$. An equivalent
definition of the symplectic fibration is that $\phi_{\alpha \beta} \subset
\Symp(M,\omega_0)$. Indeed, given such trivialization the form $\Lambda_0$
is obtained via a partition of unity from canonical forms on $\pi^{-1}
(U_{\alpha})$ such that it restricts on each fiber $M_b$ to
$\omega_b=\phi_{\alpha}(b)^* \omega_0$.

Given a symplectic fibration, we consider the {\it associated cohomological
respectively homological bundles} $H^*(M,\R) \rightarrow \Q^* \rightarrow B$
and $H_*(M,\Z) \rightarrow \Q_* \rightarrow B$. These are obtained by simply
considering the same atlas for the base and automorphisms naturally induced by
the maps $\phi_{\alpha}$ on homology respectively cohomology.

In a similar manner one constructs an associated bundle $\Jj (B,M)$ whose
fiber over each $b$ is the space of almost complex structures $J$ on $M$ that
are compatible with $\omega$. As explained in Le-Ono \cite{LO}, since the
fibers are contractible, one can always pick a section $b \rightarrow J_b$ in
this bundle.

The above alternative descriptions of a symplectic fibration imply that there
exist constant sections $s^{[\omega_0]} :B \rightarrow Q^*$ with the value $
[\omega_0] \in H^*(M,\R) $ and $s^{[c_1]} :B \rightarrow Q^*$ that takes the
integer value $c_1(M,\omega_0) \in H^2 (M, \R)$.

We say that a symplectic fibration is a {\bf hamiltonian fibration} if the
structure group further reduces to $\rm{Ham}(M,\omega_0) \subset
\Symp(M,\omega_0)$.

By a result of Guillemin
and Sternberg \cite{GS}, a symplectic fibration with a simply connected base $B$ is hamiltonian if and only if  there exist a {\it closed extension $\Lambda_0 \in \Omega^2(Q)$}. Moreover, a result of Thurston \cite{MS1} (page 333) guarantees that if the
base $B$ carries a symplectic form $\sigma_B$, then for $t$ sufficiently large
the form $\Lambda_0 $ can be chosen symplectic representing the class
$[\omega_0] + t [\pi^* \sigma_B]$.

If $\pi_1 B$ acts trivially on the associated fibration
$H_*(M,\Z) \rightarrow \Q_* \rightarrow B$(e. g. if $B $ is simply connected), then for each $D \in H_2(M,\Z)$ there also exists a constant section $s_D :B
\rightarrow Q_* $ that takes the value D.

Let us consider a symplectic deformation $(M,\omega_{\lambda})_{\lambda \geq
0}$ of the symplectic structure $(M,\omega_0)$.


\begin{definition}\label{def:fiberwisedef}
  We say that a continuous one parameter family of vertically closed 2-forms
  $(\Lambda_{\lambda})_{ \lambda \geq 0} $ on $Q$ that satisfy the conditions
  in definition (\ref{def:Fibration}) for symplectic fibers $(M,
  \omega_{\lambda})$, represents a {\it fiberwise symplectic deformation based
  on} the family $(M,\omega_{\lambda})_{\lambda \geq 0}$.
\end{definition}

These fibrations carry {\it vertical almost complex structures}
$\Tilde{J_{\lambda}}$. That is, almost complex structures on $Q$ taming
$\Lambda_{\lambda}$ and compatible with the fibration. We will refer to such
triples $(Q, \Lambda_{\lambda}, \Tilde{J_{\lambda}})$ as compatible to the
symplectic fibration with fiber $(M,\omega_{\lambda})$.

\subsection{ Some topological aspects of the symplectomorphism groups}

In the rest of the paper we will study Samelson products in the
symplectomorphism groups $G_{\lambda} = \Symp(M,\omega_{\lambda}) \cap
\Diff_0(M)$.

We will use greek letters such as $\gamma$ for $S^k$-cycles in the
symplectomorphism groups $G_{\lambda}$ and use the notation
$\widetilde{\gamma}$ for their homotopy classes in $\pi_k
G_{\lambda}$. $E(\gamma)$ and $\widetilde{E}(\gamma)$ will be the
corresponding $S^{k+1}$-cycle in the classifying space and respectively its
homotopy class in $\pi_{k+1} BG_{\lambda}$.

There is no direct inclusion of elements from $G_{\lambda}$ in
$G_{\lambda+\epsilon}$. The following proposition provides an adequate
substitute:

\begin{proposition}{Bu\c se \cite{B}}\label{prop:Kink}
  Denote by $G_{[0,\infty ) }  = \cup_{\lambda>0} G_{\lambda} \times
  \{\lambda\}\subset {\rm Diff}(M) \times [0,\infty)$.

  Consider $K$ an arbitrary compact set in $ G_{\lambda}$. Then there is an
  $\epsilon_{K}>0$ and a continuous map $h:[- \epsilon_K,\epsilon _K] \times K
  \rightarrow G_{[0,\infty )}$ such that the following diagram  commutes
\begin{equation}\label{germ}
\xymatrix{ h:[- \epsilon_K,\epsilon_K]\times K \ar[r] \ar[d]^{pr_1} &
  G_{[0,\infty)}\ar[d]^{pr_2} \\ [- \epsilon_K,\epsilon_K] \ar[r]^{incl} & (-
  \infty,\infty). }
\end{equation}
 Moreover, for any two maps $h$ and $h'$ satisfying this diagram and which coincide on $0 \times K$,
  there exists, for small enough $\epsilon' > 0$, a homotopy $H:[0,1] \times
  [- \epsilon', \epsilon'] \times K \rightarrow G_{[0,\infty]}$ between $h$
  and $h'$ which also satisfies $H \circ pr_2= pr_1 \circ incl$.
\end{proposition}

Let $\gamma_0:S^k \rightarrow G_{0}$ be a cycle in $G_{0}$. An {\bf extension} $\gamma_{\lambda}, \lambda \geq 0$ of $\gamma_0$ is a smooth family of cycles
$\gamma_{\lambda}:S^k \rightarrow G_{\lambda }$ defined for small $\lambda$
 and satisfying (\ref{germ}).

\begin{definition}
  We say that an element $\widetilde{\gamma}_0 \in \pi_*G_{0}$ is {\bf
  fragile} if it admits a null homotopic extension to the right $0=
  \widetilde{\gamma}_{\lambda}\in \pi_* G_{\lambda }$, for $\lambda > 0$. The element
  is said to be {\bf robust} if it admits an essential extension to the right
  $0 \neq \widetilde{\gamma}_0\in \pi_* G_{\lambda }$.

A continuous family $\gamma_{\lambda} :B \rightarrow G_{\lambda} , \lambda >0$
is {\bf new} if it is not the extension of a map $\gamma_0 :B \rightarrow
G_{0} $.
\end{definition}

\subsection{ Whitehead and Samelson products}

Consider a topological group $G$ and its classifying space $X= BG$ with
$\Omega X= G$. Any {\it Whitehead products} can be introduced for an arbitrary
topological space X.

Consider elements $\eta_i: S^{j_i} \rightarrow G$ representing elements
$\widetilde{\eta_i}$ in $\pi_*(G)$, and their suspensions $E(\eta_{i}) :
S^{j_i} \rightarrow B G$.

The {\it Samelson products} $ [\widetilde{\eta}_{1},\widetilde{\eta}_{2}]_s
\in \pi_{j_1+j_2}(G)$ are given by the quotient of the commutator $
[\eta_1,\eta_2]:S^1 \times S^1 \times S^{j_1+j_2} \longrightarrow G$

\begin{equation}\label{comm1}
  [\eta_1,\eta_2](s,t)=\eta_1(s)\eta_2(t)\eta_1(s)^{-1} \eta_2(t)^{-1}
\end{equation}
to $S^{j_1+j_2}= S^{j_1} \times S^{j_2} /S^{j_1} \vee S^{j_2}$.

The ordinary {\it second order Whitehead product}
$[\widetilde{E}(\eta_1),\widetilde{E}(\eta_2)]_w \in\pi_{j_1+j_2+1}(B
G_{\lambda}) $ is given by the obstruction to extending the wedge map
$E(\eta_1) \vee E(\eta_2) : S^{j_1+1} \vee S^{j_2+1} \rightarrow B G$ to a map
with the domain $S^{j_1+1} \times S^{j_2+1}$.  A classical result states that
$ [\widetilde{\eta}_1,\widetilde{\eta}_2]_s $ is, up to a sign, the
desuspension of $[\widetilde{E}(\eta_1),\widetilde{E}(\eta_2)]_w $.

  
Following \cite{P}, the $k^{\rm th}$ order higher Whitehead products
$[\widetilde{E}(\eta)_{j_1},\ldots, \widetilde{E}(\eta)_{j_k}]_w$ is a
(possibly empty) subset of homotopy classes in $ \pi_{r-1}(BG)$, with $r=j_1+1
+\ldots + j_k+1$, defined as follows:
 
Let $P= \Pi_{i=1}^k (S^{j_i+1})$. The {\it fat wedge} product $T$ is the $
r-2$ skeleton inside $P$ and consists of all the $k$-tuples $(x_1,\ldots,x_k)$
of points in $P$ such that at least one of their coordinates coincides the
coordinate of a base point $x_1^0$. Clearly $P$ is obtained from $T$ by
attaching an $r$ dimensional cell with an attaching map $a: S^{r-1}
\rightarrow T$ also called the universal Whitehead product.

Given the set of homotopy classes $\widetilde{E}(\eta_i)$ we have the
following wedge map, unique up to homotopy:
 
\begin{equation}\label{def:wedgemap}
  g=E(\eta_1) \vee\ldots \vee E(\eta_k):S^{j_1+1} \vee \ldots \vee S^{j_k+1}
  \longrightarrow BG
\end{equation}

such that $g \circ i_i= E(\eta_{i})$ with $i_i$ the obvious inclusions.  

Consider now the canonical inclusions $  i: S^{j_1+1} \vee \ldots \vee
S^{j_k+1}\rightarrow T$ and take the set of all possible extensions of $g$ 
  
\begin{equation}\label{ext}
  {\cal W}:= \{ \bar{g}| \bar{g} :T\longrightarrow BG,  \bar{g} \circ i = g\}
\end{equation}

Then the $k^{\rm th}$ order higher Whitehead product is defined as the set of
elements in $ \pi_{r-1}(BG)$ given by the maps $a \circ \bar{g} :S^{r-1}
\rightarrow BG$ for all possible extensions $\bar{g} \in {\cal W}$ and
canonical attaching maps $a:S^{r-1} \rightarrow T$. ${\cal W}$ is nonempty if
and only if all the lower Whitehead products contain the element $0$. It is
immediate that the set of elements in $[\widetilde{E}(\eta_1),\ldots,
\widetilde{E}(\eta)_k]_w$ represents the obstructions to extending all
possible maps $ \bar{g}$ to the product $P$.

 If one is interested only in those homotopy elements in Whitehead products that have infinite order those can  be obtained as Whitehead products  in a space  $X_{\emptyset}$ called the {\it rationalization} of $X$, or equivalently, localization at $\emptyset$ (cf. \cite{AA} and references therein).
There exist localization maps 
\begin{equation}
e: X \rightarrow X_{\emptyset}  
 \end{equation}
such that any for any $x \in  [e_*(\widetilde{E}(\eta_1)),\ldots,
e_*(\widetilde{E}(\eta)_k)]_w \subset \pi_N X_{\emptyset}$ there exist an integer numbers $M_i$ such that $M x= e_* z$, with $ z \in  [M_1(\widetilde{E}(\eta_1)),\ldots,
M_k(\widetilde{E}(\eta)_k)]_w$.

The Whitehead products between elements $e_*(\widetilde{E}(\eta))$  in the rationalization $X_{\emptyset}$ are called {\it rational Whitehead products}. These products are multilinear. In light of the above correspondence and since we will be interested in nontrivial elements of infinite order only up to multiplication with a factor, we will often say that the rational Whitehead products considered are of elements in $\pi_* BG \otimes \Q$. This correspondence can be well formalized if we consider other definitions of Whitehead products cf. Allday \cite{A1}, who defines rational Whitehead products on the graded differential Lie algebra  $\pi_*BG_{\lambda} \otimes \Q$. (see remark \ref{lolo}).

\begin{definition}(Andrew-Arkovitz \cite{AA}) 
  We say that the $k^{\rm th}$  order  (rational Whitehead product
  $[\widetilde{E}(\eta_1),\ldots, \widetilde{E}(\eta)_k]_w$ {\it vanishes} if
  it only contains the trivial element.
 \end{definition}

\begin{definition}
  We call $r \geq 2$ {\it the rational minimal Whitehead order} of a topological space
  $X$ if it is the minimal order in which there exists a nonvanishing rational
  Whitehead product.
\end{definition}

We will need the following result:
 
\begin{proposition}\label{arkovitz} 

     (Andrew-Arkovitz \cite{AA}) If each group homotopy group $\pi_* X_{\emptyset}$ of the rationalization of a space $X$ is finitely generated and if $r$ is the minimal rational Whitehead order
          then any rational Whitehead product of order $r$ contains exactly
          one element.
    
\end{proposition}

We use these products to introduce the higher order Samelson products. In the
present paper we will only consider the case when all elements $E(\eta_i)$ are
the same and are given as suspensions of a circle map $\gamma : S^1
\rightarrow G$. In this situation we will use, for brevity, the notation
$W^{(k)}(\widetilde{E}( \gamma))$ to denote the $k^{\rm th}$ order Whitehead
product $[\widetilde{E}(\gamma),\ldots,\widetilde{E}(\gamma)]_w \subset \pi_{2
k - 1} BG$. The $k^{\rm th}$ order Samelson product
$S^{(k)}(\widetilde{\gamma}) := [\widetilde{\gamma}, \ldots,
\widetilde{\gamma}]_s^{(k)}$ will be a set in $\pi_{2k-2}(G)$ consisting of all
the desuspensions of the elements in the set $W^{(k)}(\widetilde{E}(
\gamma))$.

\section{Whitehead products as obstructions to the existence of symplectic
         fibrations}\label{sec2}

\subsection{Towers of symplectic fibrations}
The results in this section will hold for any continuous family of symplectic structures $(M, \omega_{\lambda})$ on a compact manifold $M$ with symplectomorphism groups $G_{\lambda} = \Symp(M,\omega_{\lambda}) \cap
\Diff_0(M)$.

To each circle map $\gamma_{\lambda}: S^1 \rightarrow G_{\lambda}$ one can
associate a symplectic fibration $P_{\gamma_{\lambda}} \rightarrow S^2$, with
fiber $(M,\omega_{\lambda})$, obtained by clutching $ M\times D^+$ and $M
\times D^-$ via the identification $(x,z) =(\gamma_{\lambda}(z) (x), 1/z)$
whenever $z \in \partial D^+$.  $P_{\gamma_{\lambda}}$ is determined up to
symplectic isotopy by the homotopy class $\widetilde{\gamma}_{\lambda}$.

The main proposition of this section shows that triviality of certain
Whitehead products yields a construction of {\it towers of symplectic
fibrations} built on the fibration $P_{\gamma_{\lambda}}$, and that these
towers behave well under symplectic deformations. 

\begin{proposition}{}\label{Whitehead}

  {\renewcommand{\labelenumi}{{\bf \alph{enumi})}}
  \begin{enumerate}
    \item Fix  $\lambda$ and assume that there exist a map  
$\gamma_{\lambda}: S^1 \rightarrow G_{\lambda}$, that yield a nontrivial $\widetilde{\gamma}_{\lambda} \in \pi_1 G_{\lambda} \otimes \Q $ for which all rational  Whitehead products of orders $k $ smaller or equal that a given $p$  vanish:
 \begin{equation}
 \{0 \}= W^{(k)}(\widetilde{E}( \gamma_{\lambda})), k \leq p
\end{equation}
Then we can build on $\widetilde{\gamma}_{\lambda}$ a tower of symplectic
          fibrations of length $p$:
          \begin{equation}\label{tower}
            \xymatrix{
              (Q'^{(1)},{\Lambda'}_{\lambda}^{(1)})
                \ar@{^{(}->}[r]^i \ar[d]^{\pi _{(1)}} &
              (Q'^{(2)},{\Lambda'}_{\lambda}^{(2)})
                \ar@{--}[r] \ar[d]^{\pi _{(2)}} &
              (Q'^{(p-1)},\Lambda_{\lambda}^{(p-1)})
                \ar@{^{(}->}[r]^i \ar[d]^{\pi_{(p-1)}} &
              (Q'^{(p)}, \Lambda_{\lambda}^{(p)} )
                \ar[d]^{\pi _{(p)}} \\
              \C P^1 \ar@{^{(}->}[r]^i &
              \C P^2 \ar@{--}[r]^i &
              \C P^{p-1} \ar@{^{(}->}[r]^i & \C P^p  }
          \end{equation}
          where the forms $\Lambda_{\lambda}^{(i)}$ are vertically closed
          2-forms on $Q^{(i)}$ as in \ref{tower}, and $ Q^{(1)}$ is a clutching fibration $P_{\gamma'_{\lambda}}$, for some $\gamma'_{\lambda} $ homotopy equivalent to a power of $\gamma_{\lambda}$.

          In this diagram the morphisms preserve the fibration structure.

    \item Assume now that there exist some other tower of length $p$ as in \ref{tower} built on the element $\widetilde{\gamma}_{\lambda} \in \pi_1 G_{\lambda} \otimes Q$.
Then $  0 \in  W^{(k)}(\widetilde{E}( \gamma_{\lambda})) , k \leq p $ and hence the rational Whitehead product $ W^{(p+1)}(\widetilde{E}( \gamma_{\lambda})) $ is defined. 

Moreover, if this tower and any of its {\it N coverings}, obtained by taking an $N$ covering of $Q^{(k)}_{\lambda}$ at each step, are obstructed to extend to  towers of length $p+1$ then the rational Whitehead product  $ W^{(p+1)}(\widetilde{E}( \gamma_{\lambda})) $ must contain a nontrivial element.

\item {\bf (Extension with respect to the parameter)}   Consider now a continuous family of homotopically nontrivial circle
          maps $\gamma_{\lambda}: S^1 \rightarrow G_{\lambda} , \lambda \geq
          0$ that yield a family of nontrivial robust homotopy elements
          $\widetilde{\gamma}_{\lambda} \in \pi_1 G_{\lambda} \otimes \Q$. 

Then for  any existing tower of length $s$ as in \ref{tower} at $\lambda=0$
          must  extend
          continuously to  towers of length $s$ as in \ref{tower} built on $\widetilde{\gamma}_{\lambda}$ for small $ \lambda >0$ as in the
            following diagram
          \begin{equation}\label{towerparameter}
          \xymatrix{
            (Q^{(1)} \times [0,\epsilon_1),\Lambda^{(1)}) \ar@{^{(}--}[r]^i
              \ar[d]^{\pi _{(1)}} & (Q^{(s-1)} \times
              [0, \epsilon_{s-1}), \Lambda^{(s-1)}) \ar@{^{(}->}[r]^i
              \ar[d]^{\pi _{(s)}}&(Q^{(s)}\times [0, \epsilon_{s}), \Lambda^{(s)} )
              \ar[d]^{\pi _{(s)}} \\ 
            \C P^1 \times [0, \epsilon_{1}) \ar@{^{(}--}[r]^i & \C P^{s-1} \times
              [0, \epsilon_{s-1}) \ar@{^{(}->}[r]^i  & \C P^{s} \times [0, \epsilon_{s})}
          \end{equation}

        where $\epsilon_k >0$, and at each level $k$ the morphisms $\pi _{(k)} $ commute with
          the projections on the second factors $[0, \epsilon_k)$ and the two forms $\Lambda ^{(k)}$ restrict to the symplectic forms
          $\Lambda_{\lambda}^{(k)}$ on $Q^{(k)} \times \{\lambda \}$.
    \item {\bf (Hamiltonian case)} The tower of symplectic fibrations
          $Q^{(p)}$ is a tower of {\it hamiltonian fibrations } if and only if
          $\gamma_{\lambda} :S^1 \rightarrow \rm{Ham}(M,\omega_{\lambda})$. In
          this case the forms $\Lambda_{\lambda}^{(i)}$ can be chosen
          symplectic.
  \end{enumerate} }
\end{proposition}

\begin{proof}{ of Proposition \ref{Whitehead}}
  Let $P^{(k)} = (S^2)^k$ and $T^{(k)}$ its corresponding fat wedge. Recall
  that there is a covering map
  \begin{equation}\label{cpsym} 
    pr_{(k)} : P^{(k)} \rightarrow  \C P^k= P^{(k)} /\Ss_ k
  \end{equation}
  where $\Ss_k$ is the $k^{\rm th}$ group of permutation.

  We will denote by $h_{(k)} : S^{2k+1}\rightarrow \C P^k $ the maps used to
  attach a $(2k+2)$-dimensional cell to $\C P^k $ in order to obtain $\C
  P^{k+1}$.

  Consider the universal fibration $EG_{\lambda} \rightarrow BG_{\lambda}$ and
  let $M_{G_{\lambda}}= M \times _{G_{\lambda}} EG_{\lambda}$. Well known
  properties of classifying spaces imply that, up to symplectic isotopy, all
  symplectic fibrations with fiber $(M, \omega_{\lambda})$ and base $B$ are
  obtained as $f^*(M_{G_{\lambda}})$ for some homotopy class of maps $f :B
  \rightarrow BG_{\lambda}$.  In particular the clutching fibration
  $P_{\gamma_{\lambda}}$ is just $(E(\gamma_{\lambda}))^* (M_{G_{\lambda}}).$
  Therefore the existence of a tower of fibrations with basis $B^{(1)} \subset
  B^{(2)} \cdots \subset B^{(p)}$ is equivalent with the existence of maps
  $\phi_{\lambda,k} : B^{(k)} \rightarrow BG_{\lambda}$ that commute with the
  inclusions $i : B_k \rightarrow B_{k+j}$.

  In order to prove the direct implication in part (a) for $p > 1$, first
  assume that for an integer $p>1$ and a given value $\lambda$, $0 \in
  W^{(p)}({\widetilde E( \gamma_{\lambda})})$. Then clearly $0 \in
  W^{(k)}({\widetilde E( \gamma_{\lambda})})$ for all $k \leq p$. Therefore
  the wedge map $E( \gamma_{\lambda}) \vee E( \gamma_{\lambda})$ admits
  extensions
  \begin{equation}\label{mapwed}
    g_{\lambda,(k)} : T^{(k)}  \rightarrow BG_{\lambda} , 1 \leq k \leq p
  \end{equation}
  which commute with the inclusions $i : T^{(k)} \rightarrow T^{(k+j)}$.

  A map defined on a subset of $ P^{(k)}$ invariant under the $\Ss_k$
  action is {\it symmetric } if it commutes with the action.

  \noindent{\bf Claim:}
  {\renewcommand{\labelenumi}{(\arabic{enumi})}
  \begin{enumerate}
    \item The maps $g^{sym}_{\lambda,(k) }$ in \eqref{mapwed} can be chosen
          {\it symmetric} and they extend to symmetric maps
        \begin{equation}
          g_{\lambda,(k)}^{ext}: P^{(k)} \rightarrow BG_{\lambda}.
         \end{equation}
    \item There exist maps $f_{\lambda,(k)} : \C P^k \rightarrow BG_{\lambda}
          $ that commute with the inclusions $i : \C P^k \rightarrow \C
          P^{k+j}$, such that $f_{\lambda,(1)} =E(\gamma_{\lambda})$ and
          $g_{\lambda,(k)}^{ext} = f_{\lambda,(k)} \circ pr_{(k)}.$
  \end{enumerate} }

  We will use induction to prove the claim:

  {\bf Proof of the claim for k=2 :} Take $g_{\lambda,(2) }=
  E(\gamma_{\lambda}) \vee E(\gamma_{\lambda})$, clearly symmetric.

  The obstruction to extend the map $ f_{\lambda, (1)}=: E(\gamma_{\lambda})$
  from $S^2=\C P^1 $ to $\C P^{2}$ is given by the homotopy class $
  [E(\gamma_{\lambda}) \circ h_{(1)}] \in \pi_3 BG_{\lambda}$ and it satisfies
  $2[ E(\gamma_{\lambda})\circ h_{(1)}] = [g_{\lambda,(2) } \circ a_{(2)}]=
  W^{(2)} (\widetilde{E}(\gamma_{\lambda}))= 0$, where $a_{(2)}: S^{3}
  \rightarrow P^2$ is the universal Whitehead product map, used to attach the
  top cell of dimension $4$ on $T^{(2)}$ to obtain $P^{(2)}$.

  But since we work rationally, at the expense of replacing $\gamma_{\lambda}$ with a multiple we can kill torsion in the obstructions of the extending maps from $\C P^k $ to $\C P ^{k+1}$ and hence we conclude that $[ E(\gamma_{\lambda})\circ h_{(1)}] $ must also be zero.  Therefore we can extend the map $f_{\lambda,(1)}$ to a
  map
  $$f_{\lambda,(2)} :\C P^2  \rightarrow BG_{\lambda}$$
  Then we take the map $g_{\lambda,(2)}^{ext}: P^{(2)} \rightarrow
  BG_{\lambda}$ to be $g_{\lambda,(2)}^{ext}= f_{\lambda,(2)} \circ pr_{(2)}$,
  which is clearly symmetric and extends $g_{\lambda,(2)}$.

  {\bf Proof that the claim for k implies the claim for k+1:}

  We have that $T^{(k+1)} =\bigvee _{j=0}^{j=k} P^{(k)}_j$ where $P^{(k)}_j$
  is an identification of the product $ P^{(k)}$ with the space of
  $(k+1)$-tuples that have the coordinate in position $j$ at the base point
  $x_j$.

  By the induction step, we already have $k+1$ identical copies of the
  symmetric map $g_{\lambda,(k)}^{ext} $ and a map $f_{\lambda,(k)}$ with
  $g_{\lambda,(k)}^{ext} = f_{\lambda,(k)} \circ pr_{(k)} $ which give (using
  the relation above) a symmetric map $g_{\lambda,(k+1)}: T^{(k+1)}
  \rightarrow BG_{\lambda} $.

  Moreover, the obstruction to extend the latter map to the product is $[
  g_{\lambda,(k)} \circ a_{(k)} ] = N [ f_{\lambda,(k)} \circ h_{(k)} ] =0 $
  by hypothesis.

  Again as before we may conclude that  $[ f_{\lambda,(k)} \circ h_{(k)} ] = 0$ and hence the map $
  f_{\lambda,(k)} $ can be extended to $f_{\lambda,(k+1)}: \C P^{k +1}
  \rightarrow BG_{\lambda}$. As before, we define $g_{\lambda,(k+1)}^{ext}=
  f_{\lambda,(k+1)} \circ pr_{(k+1)}$, which is a symmetric extension of
  $g_{\lambda,(k+1)}$.  \QED

From point (2) of the claim we obtain the tower of fibrations
\eqref{tower}. Note that the forms $\Lambda_{(k)}$ can be chosen to extend one
another by defining them as pull-backs from the universal fibration.

To prove point (b) let us consider a tower of
fibrations as in (\ref{tower}) of length $p$ extending $P_{\gamma_{\lambda}}$.
This gives a sequence of maps $f_{\lambda,(k)} : \C P^k \rightarrow
BG_{\lambda} $, $1 \leq k \leq p$ commuting with the inclusions that extend
$E(\gamma_{\lambda})$. If we take $g_{\lambda,(k)}^{ext}= f_{\lambda,(k)}
\circ pr_{(k)}$ then it is immediate that these maps are just extensions to
the products $P^{(k)}$ of the recurrently constructed maps $g_{\lambda,(k)} :
T^{(k)} \rightarrow BG_{\lambda} $, $1 \leq k \leq p$, with $
g_{\lambda,(2)}=E(\gamma_{\lambda}) \vee E(\gamma_{\lambda})$. But by
definition this implies that $0 \in W^{(k)}({\widetilde E(
\gamma_{\lambda})})$ for all $k \leq p$.

To show the remaining part of point (b) let us assume that  $f_{\lambda,(p)} : \C P^p \rightarrow BG_{\lambda} $
cannot be extended over $ \C P^{p+1}$. This implies that $ [ f_{\lambda,(p)} \circ h_{(p)} ] \neq 0 $. We know that the obstruction to extend the map $g_{\lambda,(p)}$ satisfies  $[
  g_{\lambda,(p)} \circ a_{(p)} ] = M [ f_{\lambda,(p)} \circ h_{(p)} ] $.
  Again, if we work rationally we can insure (by considering a covering of the given fibration as in the hypothesis) that the homotopy class $[ f_{\lambda,(p)} \circ h_{(p)} ]$ is of infinite order and hence  $[ g_{\lambda,(p)} \circ a_{(p)} ] \neq 0$.

Point (c) is obtained by applying the Proposition \ref{prop:Kink}
for all the maps $g_{0,(k)} : T^{(k)} \rightarrow BG_0$. Note that if maps from $G_{0} \rightarrow G_{\lambda}$ exist for all $\lambda$ then the fibrations extend for all parameters $\lambda$.

For point (d) let us first observe that we can replace $G_{\lambda} $ with its
subgroup $\rm{Ham}(M,\omega_{\lambda})$ and repeat point (a) to argue that we
have a tower of hamiltonian fibrations.  Moreover, the forms $ \Lambda_{(k)}$
can be chosen symplectic by taking a choice of closed forms $\Lambda'_{(k)}$
as in point (a). Then we replace them with $ \Lambda'_{(k)} \oplus s \cdot
\omega_{\C P^k}$ for large enough $s$.
\end{proof}
\QED 
\begin{remark}\label{remark2}
{\renewcommand{\labelitemi}{$\bullet$}
\begin{itemize}
  \item We do not need the family $\gamma_{\lambda}$ to exist for all values
        $\lambda \geq 0$. One can use Proposition \ref{prop:Kink} and restate
        the results for $\lambda $ in small intervals $(0, \epsilon_{K})$.
  \item The tower of fibrations we use in Theorem \ref{Whitehead} does not
        allow us to keep track of the torsion elements in
        $\pi_*G_{\lambda}$. It would be interesting to see if a version of
        these results could be set up in order to keep track of the torsion
        elements and hence find all Samelson products in $\pi_*G_{\lambda}$.
  \item One can use different types of symplectic fibrations to decide whether
        Samelson products between distinct robust elements, not necessarily in
        $\pi_1 G_{\lambda}$, are nontrivial. G.-Y. Shi \cite{Sh} has an
        approach in this direction.
\end{itemize} }
\end{remark}

\section{Parametric Gromov-Witten invariants and Whitehead products}
Parametric Gromov-Witten invariants count vertical almost holomorphic maps in
a symplectic fibration $(M, \omega_{\lambda}) \rightarrow Q \rightarrow B$
endowed with a compatible triple $(Q,\Lambda_{\lambda}, \tilde{J})$. Le-Ono
\cite{LO} and P. Seidel \cite{S1} have used them before to detect robust
elements in the symplectomorphism groups. By contrast, we will use them to detect fragile elements.

We will carefully exploit their properties of being fiberwise symplectic
deformation invariants, in our cases of interest, to show how, combined with
Proposition \ref{Whitehead}, they yield nontrivial Whitehead products.

 Equivariant Gromov-Witten invariants are a special case of the parametric Gromov-Witten invariants that are computed on manifolds $(M, \omega_{\lambda})$ that admit hamiltonian actions by compact Lie group $H$.  We look at the case $H \approx S^1$ and treat it in a separate section.

\subsection{ Definition and properties of parametric Gromov-Witten invariants}
We will first make a summary of their defining properties. We will use results
from Li-Tian \cite{LT} as well as results from Le-Ono \cite{LO}.

Assume that the symplectic fibration $\pi:Q \rightarrow B$ with fiber
$(M,\omega)$, admits {\it a closed extension} $\Lambda$ of $\omega$. Then as
explained in Subsection \ref{sect:symplecticfibration} we may consider a
section $\Tilde{J}: B \rightarrow \Jj(B, M) $ that provides an almost complex
structure on each fiber $M_b$ compatible with the symplectic form $\omega_b$.

For $2g+m > 2$ let $\Mm _{g,m}$ be the moduli space of genus $g$ Riemann
surfaces with $k$ marked distinct points. As usual ${\overline \Mm _{g,m}}$
represents the $(3g -3 +m)$-dimensional orbifold consisting of Riemann
surfaces of genus $g$ with at most rational double points different from the
marked points; that is the Deligne-Mumford compactification of $\Mm _{g,m}$.

Fix a homology class $D \in H_2(M, \Z)$ and assume that $D$ yields a constant
section $s_D : B \rightarrow H_2(Q,\Z)$ in the corresponding homological
bundle.  This will be the case if, for instance, the base $B$ is simply
connected.

For any $(\Sigma_g, x_1 \cdots x_m) \in {\overline \Mm _{g,m}}$ the map
$f: (\Sigma_g, x_1 \cdots x_m) \rightarrow Q$ is a {\it vertical stable map}
if its image is contained in a fiber $Q_b$ and the following two conditions
are satisfied:

{\renewcommand{\labelenumi}{(\arabic{enumi})}
\begin{enumerate}
  \item Any irreducible component $\Sigma_{irred}$ of genus $0$ on which $f$
        is constant must contain at least two marked points
  \item Any irreducible component $\Sigma_{irred}$ of genus $1$ on which $f$
        is constant must contain at least one marked point
\end{enumerate} } 

Let $j \in \Teich(\Sigma)$ be an arbitrary complex structure on $\Sigma$.

Note that the condition $2g + m >2$ is not mandatory. A vertical map $f$ with
$\im (f) \subset Q_b$ is $J_b$-holomorphic if there is an arbitrary complex
structure $j \in \Teich(\Sigma)$ on $\Sigma$, such that $\pbar_{J_b} (f)
=\frac1 2 (df + J_b \circ df \circ j )=0$.

Consider $f: (\Sigma, x_1 \cdots x_m) \rightarrow Q$ and $f': (\Sigma', x'_1
\cdots x'_m) \rightarrow Q$.

$(b, f,x_1, \cdots x_m )$ is equivalent to $(b', f', x'_1, \cdots x'_m)$ if
$b=b'$, both $\im(f)$ and $\im(f')$ are contained in the same fiber $Q_b$, and
there is a biholomorphism $\phi :\Sigma \rightarrow \Sigma'$ that takes marked
points to marked points, nodal points to nodal points (and hence irreducible
components to irreducible components), and such that $f \circ \phi =f'$.

Let $ {\overline \Ff}_{g,m}^{l}(Q, \tilde{J}, s_D) $ be the moduli space
of equivalence classes of triples $[b, f, x_1,\cdots, x_m]$ as above such that
$f$ is $C^l$ smooth and $[\im(f) ] = s_D (b) \in H_2(Q_b, \Z)$.
We denote by ${\overline \Mm}_{g,m}(Q, \tilde{J}, s_D) $ the subset of $
{\overline \Ff}_{g,m}^{l}(Q, \tilde{J}, s_D) $ consisting of $J_b$-holomorphic
stable maps.

Let $\Smooth(\Sigma) \subset \Sigma$ be the set of all non-singular points of
$\Sigma$. We denote by $ \Omega^{(0,1)} (f^*TQ_b^{vert})$ the set of all
continuous sections $\xi$ in $\Hom(T\Smooth(\Sigma) , f^* TQ_b^{vert})$ that
anticommute with $j$ and $J_b$. Any such section can be continuously extended
over the nodal points of $\Sigma$. We can construct a generalized bundle $E$
over $ {\overline \Ff}_{g,k}^{l}(Q, \tilde{J}, s_D) $ with fiber
$\Omega^{(0,1)} (f^*TQ_b^{vert})$ and consider a section in $E$ given by $\Phi
=\frac1 2 (df + J_b \circ df \circ j )$. Then $\Phi^{-1}(0)$ is exactly
${\overline \Mm}_{g,m}(Q, \tilde{J}, s_D) $. The topology on the space\\
$\Hom(T\Smooth(\Sigma) , f^* TQ_b^{vert})$ will be defined as in Li-Tian
\cite{LT}.

\begin{proposition}\label{prop:fredholm}
  For  $l \geq 2$ and the section $\phi: \Ff_{g,k}^{l}(Q, \tilde{J}, s_D)
  \rightarrow E $ as above, ${\phi^{-1}(0) = \overline \Mm}_{g,m}(Q,
  \tilde{J}, s_D) $ is compact and $\phi$ gives rise to a generalized Fredholm
  orbifold bundle with a natural orientation and index $d = 2({\rm dim\,}_{\C}
  M -3)(1-g) + 2 c_1(D) + 2m + {\rm dim\,} B$.
\end{proposition}

Following as in \cite{LT}, the above  result allows one to construct a {\it virtual moduli class}\\ $[{\overline
\Mm_{g,m}(Q, \tilde{J}, s_D)} ] ^{virt} \in H_d( {\overline \Mm_{g,m}(Q,
\tilde{J}, s_D)}, \Q)$. Let us consider now the usual evaluation map \\$ev : {\overline
\Mm_{g,m}(Q, \tilde{J}, s_D)} \rightarrow Q^m $ given by $ev([b, f, x_1,
\cdots x_m]) = (f(x_1), \cdots, f(x_m))$ as well as the forgetful map $forget:
{\overline \Mm_{g,m}(Q, \tilde{J}, s_D)} \rightarrow {\overline \Mm_{g,m}}$
whose value is the stabilized domain (collapsing unstable components) of $f$.

\begin{definition}\label{pgw}

The parametric Gromov-Witten invariants are maps 

\begin{equation}
  PGW_{g,m}^{\Tilde{J}} (Q, s_D ) :[H^*(Q;\R)]^m \times H^*({\overline
  \Mm_{g,m}}, \Q) \rightarrow \Q
\end{equation}

which, for $ \alpha \in [H^*(Q;\R)]^m$ and $\beta \in {\overline \Mm_{g,m}}$
are given as :

$$
  PGW_{g,m}^{\Tilde{J}} (Q, s_D)(\beta, \alpha) = [ {\overline \Mm_{g,m}(Q,
  \tilde{J}, s_D)} ] ^{virt} \cap (forget^* \beta \cup ev ^* \alpha)
$$

These invariants are zero unless
\begin{equation}\label{pgwdim1}
  2({\rm dim\,}_{\C} M-3)(1-g) +2 c_1(D) +2m +{\rm dim\,} B= {\rm dim\,}
  \alpha + {\rm dim\,} \beta
\end{equation}
\end{definition}

Assume that the associated smooth bundle in homology is trivial. Let us focus
on the case when $\beta=0$ and $\alpha$ is the Poincar\'e dual of a product of
$k$ cycles $a_i$ that can be represented in a fiber $Q_b$ for some arbitrary
$b$.  Then the invariants count all maps $[b, f, j, x_1, \ldots, x_m]$ (with
no restrictions on the genus $g$ domain $(\Sigma, j)$), whose homology class
is $[\im(f)] =s_D \in H_2(Q_b, \Z)$ and such that $f(x_i) $ lies in $a_i$.

We define {\it the symplectic vertical taming cone } $ \Tt(\Tilde{J})$ of a
section $\Tilde{J}$ to be the space of closed 2-forms $\Lambda$ on $Q$ that
are compatible with the symplectic fibration $\pi: Q \rightarrow B $ with
fiber $(M, \omega)$ and which satisfy the taming relation $\Lambda(v,\Tilde{J}
w) >0$ for any vectors $v,w$ tangent to a fiber $Q_b$.

As in Li-Tian \cite{LT} and Le-Ono \cite{LO}, the following properties of
parametric Gromov-Witten invariants hold:

\begin{proposition}{(Properties of parametric Gromov-Witten invariants).}
\label{propofpgw}
  Consider a symplectic fibration $\pi : Q \rightarrow B$ with fiber
  $(M,\omega_0)$, with a closed extension $\Lambda_0$ of $\omega_0$ and an
  integral homology class $D \in H_2(M,Z)$.

  {\renewcommand{\labelenumi}{(\roman{enumi})}
  \begin{enumerate}
    \item The parametric Gromov-Witten invariants $PGW_{g,m}^{\Tilde{J}} (Q,
          s_D )$ are well defined and independent of the choice of the section
          of tamed vertical almost complex structure $\Tilde{J}$ with
          $\Lambda_0 \in \Tt(\Tilde{J})$.
    \item The parametric Gromov-Witten invariants $PGW_{g,m}^{\Tilde{J}} (Q,
          s_D )$ are independent of the choice of the {\it taming} closed
          extension $\Lambda$ and hence are {\bf fiberwise symplectic
          deformation invariants } as long as the deformation is within some
          symplectic taming cone $ \Tt(\Tilde{J})$.
    \item(Le-Ono \cite{LO}) {\bf Symplectic sum formula} Let $Q=Q_1 \# Q_2$ be a fiber connected
          sum of two fibrations. Then
          \begin{equation}
            PGW_{g,0}^{\Tilde{J}} (Q, s_D )= PGW_{g,0}^{\Tilde{J_1}} (Q_1, s_D
            ) + PGW_{g,0}^{\Tilde{J_2}} (Q_2, s_D )
          \end{equation}
    \item (Le-Ono \cite{LO}) If $f: B' \rightarrow B$ is a $N$ covering map then
          $PGW_{g,m}^{\Tilde{J}} (Q, s_D ) = N \cdot PGW_{g,m}^{f^*\Tilde{J}}
          (f^*Q, s'_D ) $
  \end{enumerate} } 
\end{proposition}

\subsection{Equivariant Gromov-Witten invariants}\label{secegw} 

Equivariant Gromov-Witten invariants can be defined for any hamiltonian action
of a compact Lie group $H$ on a symplectic manifold $(M,\omega)$.  We will
restrict ourselves to the case of hamiltonian circle actions.

Consider the universal symplectic fibration $M_{S^1}= M \times_{S^1} E S^1$
with fiber $(M,\omega)$. $M_{S^1}$ consists of an infinite tower of
hamiltonian fibrations $\pi_{(k)} :M_{S^1}^{(k)}= M \times_{S^1} S^{2k+1}
\rightarrow \C P^k $. Note that $M \times_{S^1} S^{2k+1}$ comes equipped with a natural
$S^1$-invariant almost complex structure $J^{(k)}$ compatible with the
fibration that makes the map $\pi_{(k)} $ almost holomorphic, as well as with
the closed extensions $\Lambda^{(k)}$ consisting of symplectic $S^1$ invariant
forms.

We say that $M_{S^1}$ admits the vertical almost complex structure $\Tilde{J}$
if $\Tilde{J}$ restricts to an usual vertical almost complex structure on each
$M_{S^1}^{(k)}$. Similarly, we say that $\Lambda$ is a closed 2-form on
$M_{S^1}$ if it restricts to a closed 2-form on each $M_{S^1}^{(k)}$. Both
$\Tilde{J}$ and $\Lambda$ can be chosen $S^1$-invariant and compatible.

We get an $S^1$-action on the spaces of maps $ {\overline
\Ff}_{g,m}^{l}(M_{S^1}, \tilde{J}, s_D) $ and $ {\overline
\Ff}_{g,m}^{l}(M_{S^1}, \tilde{J}, s_D) $ and we can construct an {\it
equivariant virtual class} 
$$[{\overline \Mm_{g,m}(M, \tilde{J}, s_D)} ]
^{virt}_{equiv} \in H_*( {\overline \Mm_{g,m}(M_{S^1}, \tilde{J}, s_D)}, \Q)$$
 Then the equivariant Gromov-Witten invariants are maps
\begin{equation}\label{egwdef}
  EGW_{g,m}^{\Tilde{J}} (M, s_D ) :[H^*(M_{S^1}\Q)]^m \times H^*({\overline
  \Mm_{g,m}}, \Q) \rightarrow H^*(BS^1,\Q)
\end{equation}
which, for $ \alpha \in [H^*(M_{S^1},\Q)]^k$ and $\beta \in {\overline
\Mm_{g,k}}$ are given as:
$$
  EGW_{g,k}^{\Tilde{J}} (M, s_D)(\beta, \alpha) = [ {\overline \Mm_{g,k}(Q,
  \tilde{J}, s_D)} ] ^{virt}_{equiv} \cap (forget^* \beta \cup ev ^* \alpha),
$$
where the ``$\cup$'' is obtained from equivariant integration.
The following proposition gives properties of equivariant Gromov-Witten invariants:
\begin{proposition}
\label{propertiesegw}
{\renewcommand{\labelenumi}{(\roman{enumi})}}
\begin{enumerate}
  \item For any vertical $S^1$-invariant almost complex structure $\Tilde{J}$
  compatible with the fibration and for any $S^1$-invariant taming form
  $\Lambda$ on $M_{S^1},$ the invariants $EGW_{g,m}^{\Tilde{J}} (M, s_D )$ are well defined
  and independent of the choice of the invariant taming vertical almost
  complex structure $\Tilde{J}$.
 \item If the equivariant class $ \alpha = \bigoplus_{k=1}^{k=\infty}
\alpha^{(k)} \u^k \in H^*(BS^1, \Q) $ then we immediately have
\begin{equation}\label{egwinsum}
EGW_{g,k}^{\Tilde{J}} (M, s_D)(\beta, \alpha) = \bigoplus_{k=1}^{k=\infty}
EGW_{g,k}^{(k),J^{(k)}} (Q^{(k)}, s_D)(\beta, \alpha ^{(k)}) \u^k 
\end{equation}           
where $EGW_{g,k}^{(k),J^{(k)}} (Q^{(k)}, s_D)(\beta, \alpha ^{(k)}) $ are the
parametric Gromov-Witten invariants of the fibration $Q^{(k)}$ and are zero
unless
\begin{equation}\label{pgwdim2}
  2({\rm dim\,}_{\C} M-3)(1-g) +2 c_1(D) +2k +2m= {\rm dim\,} \alpha^{(k)} +
  {\rm dim\,} \beta
\end{equation}
\end{enumerate}
\end{proposition}

\subsection{ Parametric Gromov-Witten invariants and Whitehead products}
\begin{lemma}\label{GWobstruction}
  Consider a symplectic deformation $(M,\omega_{\lambda})_{\lambda \geq 0)}$
  and a homology class $D \in H_2(M,\Z) $ with $[\omega_0](D)=0$.  Assume that
  there exists a smooth symplectic fibration $\pi :Q \rightarrow B$ endowed
  with a continuous family of closed two extensions
  $(\Lambda_{\lambda})_{\lambda >0}$ of the symplectic fibers
  $(M,\omega_{\lambda})$. Moreover, assume that the maps
  $PGW_{g,m}^{\Tilde{J}} (Q_{\lambda}, s_D) $ are nontrivial.

  Then the family $(\Lambda_{\lambda})_{\lambda >0}$ {\bf cannot} extend to a
  fiberwise symplectic deformation $(\Lambda_{\lambda})_{ \lambda \geq 0} $
  based on the given family $(M,\omega_{\lambda})_{\lambda \geq 0}$.
\end{lemma}

The proof is immediate. Indeed, the existence of a $J $-holomorphic curve in a
class $D \in H_2(M, \Z)$ implies that any taming symplectic form $\omega_0$
must satisfy $[\omega_0](D)>0$. The result is a consequence of point (ii) in
Proposition \ref{propofpgw}.

We will effectively use the lemma above to show that extensions with respect
to the parameter as in Proposition \ref{Whitehead} cannot exist in the
presence of certain nontrivial PGW invariants:

\begin{corollary}\label{pgwobstructfibration}
  Assume that we are in the conditions of Proposition \ref{Whitehead} point (c) and we have a tower at level $\lambda=0$ of length $p \geq 1$. Then the resulting fibrations
  $Q_{\lambda}^{(k)}, k \leq p$, $(\lambda >0)$ obtained  by extending with respect to the parameter the fibrations $ Q_{0}^{(k)},k\leq p$  must
  satisfy $0= PGW_{g,m}^{\Tilde{J_{\lambda}}} (Q^{(p+1}_{\lambda}, s_D) $
  whenever $[\omega_0] (D)=0$.
\end{corollary}
The crux of the argument that $\{ 0 \} \neq W^{(p+1)}({\widetilde E(
\gamma_{0})})$, will be to show that some fibration $Q^{(p+1)}_{\lambda}$
obtained by extending $Q_0^{(p+1)}$ with respect to the parameter, must have a
nontrivial parametric Gromov-Witten invariant as above which would contradict
the above corollary.

\section{Ruled surfaces}\label{sec5}

A ruled surface $M_{\lambda}^g$ is the total space of the topologically
trivial symplectic fibration $(S^2 \times \Sigma_g, \sigma_{S^2} \oplus
\lambda \sigma_{\Sigma_g}) \rightarrow (\Sigma_g,
\sigma_{\Sigma_g})$. Accordingly, we let the symplectomorphism groups
$G_{\lambda}^g$ be $\Symp(S^2 \times \Sigma_g, \sigma_{S^2} \oplus \lambda
\sigma_{\Sigma_g})\cap \Diff_0(M)$.

\subsection{ Prior results}

We will present here results that are essentially contained in McDuff \cite{M}.
Let us denote by $\Ss_{\lambda} $ the space of symplectic forms that are
strongly isotopic with $\omega_{\lambda}$, and by $\Aa_{\lambda}$ the space of
almost complex structures that are tamed by some form in $\Ss_{\lambda}$.
Then there exists a fibration $G_{\lambda} \rightarrow \Diff_0 (M)
\rightarrow \Ss_{\lambda} $ and, since $\Ss_{\lambda} $ is homotopy
equivalent with $\Aa_{\lambda} $, there is also a homotopy fibration
\begin{equation}\label{acsfibration}
  G_{\lambda} \rightarrow \Diff_0 (M) \rightarrow \Aa_{\lambda}.
\end{equation}
Let $D_k = A -k F \in H_2(M^g_{\lambda},\Z)$ where $A$ and $F$ are the
homology classes of the base and the fiber respectively. The subsets
$\Aa^g_{\lambda, k} $ of $\Aa^g_{\lambda} $ consisting of almost complex
structures that admit $J$-holomorphic curves in the class $D_k$ provide a
stratification of $\Aa^g_{\lambda} $ as in the following:

\begin{proposition}{(McDuff\cite{M}) }\label{mcduffacs} 
  {\renewcommand{\labelenumi}{{\bf (\roman{enumi})}}
  \begin{enumerate}
    \item $\Aa^g_{\lambda} \subset \Aa^g_{\lambda +\epsilon} $ and hence, via
          (\ref{acsfibration}) one obtains maps $h_{\lambda, \lambda+\epsilon}
          :G_{\lambda} \rightarrow G_{\lambda+ \epsilon} $.
    \item $\Aa^g_{\lambda, k} $ is a Frechet suborbifold of $\Aa^g_{\lambda} $
          of codimension $4k-2+2g.$
    \item $\Aa^0_{\lambda} $ is constant on all the intervals $(\ell, \ell+1]$
          and $\Aa^0_{k+\epsilon}  \setminus \Aa^0_{k} = \Aa^0_{k+\epsilon, k} $.
    \item The homotopy type of $G_{\lambda}^0$ is constant for $k < \lambda
          \leq k+1$, with $k$ an integer greater than zero. For this range of
          $\lambda$ there exists a nontrivial fragile element $w_k \in
          \pi_{4k}( G_{\lambda}^0) \otimes \Q$ that disappears when $\lambda$
          passes the critical value $k+1$, while a new fragile element
          $w_{k+1} $ appears.

    \item There exists a fragile element $\rho \in \pi_2 G_{1}^1 $ that
          disappears in $ \pi_2 G_{1+\epsilon}^1. $
  \end{enumerate} }
\end{proposition}

Moreover the inclusions $i: G_{\lambda}^g \rightarrow \Diff_0(M^g) $ lift to
maps ${\tilde i} : G_{\lambda}^g \rightarrow \Dd_0^g$ where $\Dd_0^g$ is the
subgroup of diffeomorphisms that preserve the $S^2$ fibers. The following
proposition shows that all essential elements in $\pi_* (\Dd_0^g) $ are
retained in the homotopy groups of symplectomorphism groups:

\begin{proposition} {McDuff\cite{M}}\label{mcduffdiff} 
  {\renewcommand{\labelenumi}{{\bf (\roman{enumi})}}
  \begin{enumerate}
    \item The vector space $\pi_i( \Dd_0^g) \otimes \Q$ has dimension 1 when
          $i=0,1,3$ except in the cases $i=g=1$ when the dimension is 3, and
          $g=0, i=3$ when the dimension is 2. It has dimension 2g when $i=2$
          and is zero otherwise.
    \item There exist maps $\tilde{i}: G_{\lambda}^g \longrightarrow
          \Dd_0^g$ that induce a surjection on all rational homotopy groups
          for all $g>0$ and $\lambda \geq 0$. The map is actually an
          isomorphism on $\pi_i$, $i = 1, \ldots, 2g-1$ when we restrict to
          the range $\lambda > k$ where $g=2k $ or $g=2k +1$ depending on the
          parity.
    \item The map $\tilde{i}$ also gives an isomorphism on $\pi_i$ for
          $g=1, i=2,3,4,5 $ and $\lambda > 3/2$.
    \item The homotopy limit $G^g_{\infty}=\lim_{\lambda  \rightarrow \infty} G_{\lambda}^g \approx {\cal D}_0^g$ 
\end{enumerate} }

\end{proposition}

\subsection{Hamiltonian circle actions on ruled surfaces, robust elements
            and  equivariant Gromov-Witten invariants}

We will first describe all possible hamiltonian circle actions on the
manifolds $M_{\lambda}^g$. This is provided for instance in M. Audin
\cite{Au}. For these actions we will give a complete description of the
equivariant Gromov-Witten invariants that count isolated curves of genus $g$.
We also describe families of robust elements that satisfy hypothesis $H_1$
which, combined with the non-trivial count of EGW yields nontrivial Whitehead
products.

  The Lie groups $H_k \approx S^1$ act on the manifolds $M^g_{\lambda} $,
$\lambda>k$ as follows:

We denote by $\Oo(-2k)_g$ to be a holomorphic line bundle of degree $-2 k$
over the surface $\Sigma_g$, and consider the projectivized line bundles
$\pi:P( \Oo(-2k)_{g} \oplus \Oo_g) \longrightarrow \Sigma_g$. The K\"ahler
manifolds $ P( \Oo(-2k)_{g} \oplus \Oo_g)$ are endowed with naturally
integrable almost complex structures denoted by $J^{(k),g}$. Topologically,
they are just $\Sigma_g \times S^2 $ and it is easy to see that these bundles
admit a {\it holomorphic} circle action that rotates the fibers while fixing
the zero section and the section at infinity that represent the classes $A-kF$
and $A+kF$ respectively:

\begin{equation}\label{actionruled}
  \gamma_k^g: S^1 \times P( \Oo(-2k)_{g} \oplus \Oo_g) \rightarrow P(
  \Oo(-2k)_{g} \oplus \Oo_g)
\end{equation}

In coordinates, this action is given by $ e^{i t} \cdot (b,
[v_1:v_2])= (b, [e^{it} v_1: v_2])$ .

We will view the $P( \Oo(-2k)_{g} \oplus \Oo_g)$ as the symplectic manifolds
$M^g_{\lambda}$ endowed with the $S^1$-invariant taming complex structures
$J^{(k),g}$ whenever $\lambda >k$.

The circle actions \eqref{actionruled} become {\it hamiltonian} with respect
to the symplectic forms $\omega_{\lambda}$, whenever $\lambda >k$; this is for
example explained in Audin \cite{Au}. The ruled surfaces $M_{\lambda}^g$, for
$\lambda >0$, can be constructed via symplectic reduction from disk bundles $
D_a( \Oo(-2k)_{g} \oplus \Oo_g)$ with appropriate radii $a$.

This construction fails to work if $g>0$ and $\lambda=1$. In fact, it follows
from Karshon \cite{K} that the symplectic manifolds $M_{0}^g$ do not admit any
such hamiltonian $S^1$-action. Moreover, it is clear that the actions
\eqref{actionruled} cease to be symplectic whenever $\lambda \leq k$.  To
stress this distinction we will use the following notation for the
hamiltonian actions
\begin{equation}\label{actionruledsymp}  
\gamma_{k,\lambda} ^g: S^1 \times M_{\lambda}^g \rightarrow
  M_{\lambda}^g , \lambda > k
\end{equation}
 or, equivalently, 
\begin{equation}\label{actionruledsympc}
  E(\gamma_{k,\lambda} ^g): S^2 \rightarrow BH_k \subset BG_{\lambda}^g,
  \lambda > k.
\end{equation}
From Proposition \ref{mcduffdiff} we see that the cycles
$\tilde{i}(\gamma_{k,\lambda} ^g)$ are essential in $\Dd_0^g$ and
represent an element ${ \widetilde \gamma^g_k} \in \Dd_0^g \otimes Q$.

In fact, a smooth representative for ${ \widetilde \gamma^g_k} \in \Dd_0^g
\otimes Q$ can be given as $\gamma'^g_k :S^1 \longrightarrow \Dd_0^g$
\begin{equation}
  \gamma'^g_k (\theta)(w,z)=(z, \rho(R_{\theta}^z(w))
\end{equation}
where $R_{\theta}^w(z)$ rotates the fiber sphere in $S^2 \times S^2$ with an
angle $\theta$ about a point $z$ in the base sphere, and $\rho : \Sigma_g
\rightarrow S^2 $ is a covering map of degree $k$.

In the case $g=0$ the hamiltonian $S^1$-actions \eqref{actionruled} are in
fact induced from a $T^2$ {\it toric} action. $M_{\lambda}^0, \lambda \geq 1 $
can be obtained through symplectic reduction in $\lfloor \lambda \rfloor$
different ways as $M^0_{\lambda} = \C^4 // T^2$, for any $0 \leq k < \lambda$
where the two generators $\xi_1, \xi_2$ of $ T^2$ act on $\C^4$ with
weights $(1,1, 0,0)$ and $(2k,0,1,1)$. The group of toric automorphisms is a
subgroup of the symplectomorphism group and it contains a Lie subgroup $ K_k =
S^1 \times SO(3) $ for $k>0$ and $ K_0 = SO(3) \times SO(3) $ such that the
map $\pi_1 K_k \rightarrow G^0_{\lambda}$ induces an injection on the homotopy
groups. In this case we have $H_k \subset K_k$.
It follows that $\pi_*K_k$ contains the generators
$\widetilde{\gamma}^0_{\lambda, k} \in \pi_1 G_{\lambda}^0$ and
$\widetilde{\alpha}_k \in \pi_3 G_{\lambda}^0$ whenever $\lambda> k > 1$, and
$\widetilde{\alpha}$ and $\widetilde{\eta}$ in $\pi_3 G_{\lambda}^0$ whenever
$\lambda >1$.
 
\begin{lemma}\label{actionmultiple}
  Consider (1) $g>0$ and  $k \geq 1$ or (2)  $ g=0$ and $k \geq 2$.
  {\renewcommand{\labelenumi}{(\roman{enumi})}
  \begin{enumerate} 
  \item \begin{equation}
\widetilde{\gamma}^g_k = k \widetilde{\gamma}^g_1 \in \pi_1 {\cal D}^g_0 \otimes \Q
\end{equation}
\item If, in addition, we assume $\lambda >k>[ g/2] $ then the same relation takes place in $\pi_1 G^g_{\lambda} \otimes \Q$:
\begin{equation}
\widetilde{\gamma}^g_{\lambda,k} = k  \widetilde{\gamma}^g_{\lambda,1} \in \pi_1 G_{\lambda}^g \otimes \Q.
\end{equation}
\item(for $g=0$)  \begin{equation}\label{relationpi3}
            \widetilde{\alpha}_k =\widetilde{\alpha} + k^2 \widetilde{\eta}
            \in H_3(G_{\lambda},\Q)
          \end{equation}
\item There exist a
  continuous family of robust elements of infinite order
  $\gamma'^g_{\lambda, k} : S^1\rightarrow G_{\lambda} ^g $, for $\lambda
  \geq k$ which for $\lambda >k$ is homotopy equivalent with the circle maps $\gamma^g_{\lambda,k}$  given by the group action $H_k$. 

 Moreover, with the exception of the case $g=k=1$,    at the critical values $\lambda=k$ we can deduce that there are integers $M$ so that 
\begin{equation}\label{fish} 
\widetilde{\gamma}'^g_k = M \widetilde{\gamma}^g_1 \in \pi_1G^g_k \otimes \Q\end{equation}
  \end{enumerate} }
\end{lemma}

\begin{proof}{}
  The proof of (i) is an immediate adaptation of  Lemma 2.10 proved in Abreu-McDuff \cite{AM} for the case $g=0$. In fact they actually compute the difference between the two terms as a 2-torsion element.
Similarly, when we restrict to the given range for $\lambda$ the morphisms $\tilde{i}$ give an isomorphism on $\pi_1$ and hence the relation in (i) continues to hold in $\pi_1G^g_{\lambda} \otimes \Q$.
Part (iii) is also contained in Lemma 2.10 proved in Abreu-McDuff \cite{AM}.

The existence of the robust family in part (iv) is an immediate consequence of Proposition \ref{mcduffdiff}. Indeed,
  since the maps $\widetilde{i}$ induce a surjection on the first rational homotopy groups for $\lambda$ in that range and the family can be obtained by pulling back the smooth representative $\gamma^g_k$ to the symplectomorphism groups. Finally, the relation \ref{fish} follows for instance from the fact that the vector space $\pi_1 G^g_k \otimes \Q$ is one dimensional (proposition \ref{mcduffdiff} point (i)).

\end{proof}
\QED

\subsection{Equivariant Gromov-Witten invariants}

\begin{proof}{ of Theorem \ref{egwruled}}
Denote by
\begin{equation}\label{fibrationnot}
  (Q_{\lambda}^{(k),(p), g}, J^{(k),(p),g}) =
  M^g_{\lambda} \times _{H_k} S^{2p+1}
\end{equation}
the associated symplectic fibration with fiber $(M^g_{\lambda},
\omega_{\lambda} )$ endowed with the $S^1$-invariant symplectic form
$\Lambda_{\lambda}^{(k),(p)}$ and compatible almost complex structure
$J^{(k),(p),g}$.  Then according with Proposition \ref{propertiesegw} we need
to show that $EGW_{g,0}^{(p),J^{(k),(p),g}} (Q_{\lambda}^{(p),g}, s_{D_k})$ is
$\pm 1$ if $p=2k +g-1$ and zero otherwise.

The dimension condition in \ref{pgwdim2} translates into saying that
\begin{multline}
({\rm dim\,}_{\C} M_{\lambda}^g-3) (1-g) + c_1(A-kF) +2p = g-1 + (A-kF)^2
+2-2g+2p=\\ = g-1 -2k +2-2g= -2k -g +1 +2p
\end{multline}

 must be 0. Therefore all such
invariants are zero unless $p= 2k+g-1$.

In this situation there exists exactly one embedded vertical
$J^{(k),(p),g}$-holomorphic map representing $s_{D_k}$ in each fiber
$Q^{(k),(p),g} _b$ for each $b \in \C P^{p}$. More precisely, each fiber is
biholomorphic to $ P(\Oo (-2k)_g \oplus \Oo_g)$. The only possible bubbling
for vertical almost holomorphic curves in $Q^{(k),(p),g} _b$ must take place
within a fiber. It immediately follows that the only $ J^{(k),(p),g}$ maps in
each fiber representing the class $D_k$ is the zero section of the bundle $
P(\Oo (-2k)_g \oplus \Oo_g)$.  Therefore the moduli space
$\Mm_{g,0}(Q^{(k),(p),g} ,J^{(k),(p),g}, s_{D_k})$ is naturally diffeomorphic
with $\C P^{p}$.

Given such $J^{(k),(p),g} $ holomorphic map $f: (\Sigma_g, j_g) \rightarrow
(M_{\lambda} ^g, J^{(k),(p),g})$ in the class $D_k$ the linearized operator
$D \phi$ of index zero is
\begin{equation}\label{ruledoperator}
  D\phi([b,f,j_g]):T_b \C P^{p} \times C^{\infty} (f^*TM_{\lambda}^g) \times
  T_{j_g} \Teich_g \rightarrow \Omega^{(0,1)} (f^*TM_{\lambda}^g)
\end{equation}
where the component corresponding to the Teichm\"uller space appears when $g >
0$.

The actual dimension of $\Mm_{g,0}(Q^{(k),(p),g},J^{(k),(p),g}, s_{D_k})$ is
larger than its formal dimension $0$. This is because the fiberwise almost
complex structure $J^{(k),(p),g}$ is not $D_k$-{\it regular}, or equivalently,
the linearized operator (\ref{ruledoperator}) is not onto. The computation of
the invariants then follows from the following:

\begin{lemma}\label{obstructionbundle}
  {\renewcommand{\labelenumi}{(\roman{enumi})}
  \begin{enumerate}
    \item $EGW_{g,0}^{(p),J^{(k),(p),g}} (Q^{(k),(p),g}, s_{D_k}) = e({ \cal
          O}^g)$ where $ e({ \cal O}^g)$ represents the Euler class of the
          obstruction bundle $ { \cal O}^g \rightarrow
          \Mm_{g,0}(Q^{(k),(p),g},J^{(k),(p),g}, s_{D_k})$ induced by the
          section $\phi$ whose fiber over a point $[b,f,j_g]$ is given by
          $\coker D\phi([b,f,j_g])$.
    \item Whenever $p=2k+g-1$ the obstruction bundle $ { \cal O}^g \rightarrow
          \Mm_{g,0}(Q^{(k),(p),g},J^{(k),(p),g}, s_{D_k}) $ is isomorphic to
          $\Oo_{\C P^p}(-1)^p \rightarrow \C P^{p}$.
  \end{enumerate} }
\end{lemma}

\begin{proof}{}
  (i) This follows immediately from the setup in the general theory as in
  Li-Tian \cite{LT}, since in this particular case the moduli space
  $\phi^{-1}(0)$ is smooth and hence the generalized Fredholm orbifold is in
  fact a smooth vector bundle over $\C P^p$.

  (ii) Since $f$ represents the zero section in the fiber $Q_b= P( {\cal
  O}(-2k)_{g} \oplus \Oo_g)$, the vertical tangent bundle
  $T_b^{vert}(Q_b^{(k)(p),g})_{|\im f} = T (M_{\lambda}^g)_{| \im f}$ splits
  holomorphically in the direct sum $T \Sigma_g \oplus \nu_g$, where $\nu^g_k$
  is the normal bundle to the image $\Sigma_g$ of the zero section $f$. It is
  immediate that the normal bundle is in fact $ \Oo(-2k)_{g} \rightarrow
  \Sigma_g$. The operator (\ref{ruledoperator}) becomes:

  \begin{multline*}
    D\phi([b,f,j_g]): T_b \C P^p\oplus C^{\infty} (\Sigma_g,\nu^g_k ) \oplus
    C^{\infty} (\Sigma_g,T\Sigma_g) \oplus T_{j_g} \Teich_g \rightarrow \\
    \rightarrow \Omega^{(0,1)} (\Sigma_g,\nu^g_k ) \oplus
    \Omega^{(0,1)}(\Sigma_g,T\Sigma_g)
  \end{multline*}
  and hence
  \begin{multline*}
    D\phi [b,f,j_g]):T_b \C P^p\oplus C^{\infty} (\Sigma_g, \Oo(-2k)_{g})
    \oplus C^{\infty} (\Sigma_g,T\Sigma_g) \oplus T_{j_g} \Teich_g \rightarrow
    \\ \rightarrow \Omega^{(0,1)} (\Sigma_g, \Oo(-2k)_{g} ) \oplus
    \Omega^{(0,1)}(\Sigma_g,T\Sigma_g)
  \end{multline*}

  We will study the cokernel in the case $g=0$ separately. If $g>0$, then the
  component of $D\phi [b,f,j_g])$ that is not onto is
 
  \begin{equation}\label{cokernel}
    D \phi _{restr}([b,f,j_g]): C^{\infty} (\Sigma_g, \Oo(-2k)_{g})
    \rightarrow \Omega^{(0,1)} (\Sigma_g, \Oo(-2k)_{g} )
  \end{equation}
  whose cokernel is $H^{(0,1)} (\Sigma_g, \Oo(-2k)_{g})$. If we denote by
  $K_g$ the degree $2g -2$ canonical bundle over $\Sigma_g$ then, by Serre
  duality, $\coker D\phi [b,f,j_g])$ will be precisely the space of
  holomorphic sections $(H^0(\Sigma_g, \Oo(-2k)_{g}^* \otimes K_g))^*$. By the
  Riemann-Roch theorem this space has complex dimension $ 2k +2g-2 -g+1=
  2k+g-1 $.  To find out how these fibers fit together topologically in the
  obstruction bundle we need to understand what is the induced $S^1$-action on
  $(H^0(\Sigma_g, {\cal O}(-2k)_{g}^* \otimes K_g))^*$ such that
  $$
    \Oo^g = (H^0(\Sigma_g, \Oo(-2k)_{g}^* \otimes K_g))^* \times
    _{S^1} S^{2p+1}.
  $$
  Since $S^1$ acts with weight $1$ on the normal bundle $\mu_g =
  \Oo(-2k)_{g}$, and correspondingly on its dual $\Oo(-2k)_{g}^*$, the space
  of sections inherits a diagonal $S^1$-action with equal weights given by
  either $1$ or $-1$. Since it will be enough to determine the EGW up to a
  sign, we will assume for simplicity that the weights are equal to 1. Since
  $\im(f)$ is a fixed set of the canonical bundle $S^1$-action,
  $T(Q^{(k),(p),g})_{| \im f}$ is also fixed by the induced $S^1$-action and
  therefore so is the $K_g$. Hence the action on $(H^0(\Sigma_g,
  \Oo(-2k)_{g}^* \otimes K_g))^*$ is induced by the $S^1$-action with weights
  $(1, \ldots, 1)$ on $\Oo(-2k)_{g}^*$ and hence it is diagonal with weights
  $(1, \ldots, 1)$. It immediately follows that $ (H^0(\Sigma_g,
  \Oo(-2k)_{g}^* \otimes K_g))^* \times _{S^1} S^{2p+1}$ is given by ${\cal
  O}_{\C P^p}(-1)^p \rightarrow \C P^{p}$.

  In the case $g=0$ the moduli spaces involved in the computation must be of
  {\it unparameterized curves}, which means we have to quotient out the
  6-dimensional group $PGL(2,\C)$ representing the reparametrizations of the
  domain.  The linearized operator will be

  $$
    D\phi([b,f,j_0]): T_b \C P^p\oplus C^{\infty} (\Sigma_0,\nu^0_k ) \oplus
    C^{\infty} (S^2,T S^2)  \rightarrow
    \Omega^{(0,1)} (S^2,\nu^g_k ) \oplus \Omega^{(0,1)}(S^2,TS^2)
  $$
  with the cokernel given by:

  $$
    D \phi _{restr}([b,f,j]): C^{\infty} (S^2, \Oo(-2k))
    \rightarrow \Omega^{(0,1)} (S^2, \Oo(-2k) )
  $$

  A similar line of thought as above then applies. In this case the canonical
  bundle is of negative degree $\Oo(-2)$ and the fiber of the obstruction
  bundle is $(H^0(S^2, \Oo(-2k)^* \otimes \Oo(-2))^*= (H^0(S^2,\Oo(2k -2))^*$ of
  complex dimension $2k-1$.

  Hence whenever $p=2k+g-1$ we have $EGW_{g,0}^{(p),J^{(k),(p),g}}
  (Q^{(k),(p),g}, s_{D_k}) = e(\Oo^g)=$\\ $ c_n ( \Oo_{\C P^p} (-1)^p) = (c_1
  ( {\cal O}_{\C P^p} (-1))^p = 1$.
\end{proof}\QED
\end{proof}
\begin{remark}
As in Proposition \ref{Whitehead} point (b) we also need to consider towers of fibrations that are finite covers of the original ones. Note that any convering of $ \Q_{\lambda}^{(k),(p), g}$ must also have nontrivial PGW cf. Proposition \ref{propofpgw}(iii).
\end{remark}

\subsection{A non-trivial Whitehead product in the symplectomorphism group of
            {\boldmath $T^2 \times S^2$}}

\begin{proof} { of Proposition \ref{fragileruledhg}(i)}
The result will follow from:
 
\noindent{\bf Claim: } The Whitehead product
  $[\widetilde{E}(\gamma'^1_{1,1}),\widetilde{E}(\gamma'^1_{1,1})]_w \in \pi_3
  BG_{1}^1 $ is nontrivial and yields, by desuspension, a nontrivial fragile
  element $ w= [\widetilde{\gamma}'^1_{1,1},\widetilde{\gamma}'^1_{1,1}]_s \in
  \pi_2 G_{1}^1$.

Let us assume that the claim is false, therefore
$[\widetilde{E}(\gamma'^1_{1,1}),\widetilde{E}(\gamma'^1_{1,1})]_w =0 $.
 
 
Then, according with Proposition \ref{Whitehead}(iii) there exists a
continuous family of fibrations $ ( Q_{\lambda}^{(2)},\Lambda_{\lambda})
\rightarrow \C P^2 $, for $\lambda \geq 1$ sufficiently close to 1 which fits
in the tower \eqref{towerparameter}. From corollary \ref{pgwobstructfibration}
we must have

\begin{equation}\label{pgwzero}
  PGW _{1,0}^{\tilde{J} }(Q_{\lambda} ^{(2)}, s_{A-F})=0
\end{equation}

Then under the triviality assumption we claim the following:

\begin{lemma}\label{fibersumtorus}
  Take $ E^{(2)} _{\lambda} \rightarrow \C P^2 $ be a fibration obtained as a
  N-covering of the fibration $Q^{(1),(2),1}_{\lambda}$ defined in
  \eqref{fibrationnot}.
  {\renewcommand{\labelenumi}{(\roman{enumi})}}
  \begin{enumerate} 
    \item For $\lambda >1$ sufficiently close to 1, there exists a continuous
          family of fibrations $R_{\lambda}\rightarrow S^4 $ given by a family
          of elements $ \eta_{\lambda} :S^3 \rightarrow G_{\lambda}$, such
          that $E^{(2)} _{\lambda} \# R_{\lambda}$ is symplectically isotopic
          to $Q_{\lambda} ^{(2)}$.
    \item $ PGW _{1,0}^{\tilde{J} }(R_{\lambda}, s_{A-F})\neq 0$
    \item The family $ \eta_{\lambda} :S^3 \rightarrow G_{\lambda}$ is {\bf
          new}.
  \end{enumerate}
\end{lemma}

\begin{proof}{ of the lemma}
  Let us remind the reader that in order to obtain the family of fibrations
  $Q_{\lambda}^{(2)} $ we take the symmetric wedge map $g'_{1,(2)} $, its
  extensions to the product $g'^{ext}_{1,(2)} $ and $f'_{1,(2)} $ to the
  product $S^2 \times S^2$ and to $ \C P^2$ respectively and extend them
  continuously with respect to the parameter. Thus we obtain maps
  $g'_{\lambda,(2)}$, $g'^{ext}_{\lambda,(2)}$ and $f'_{\lambda,(2)}$ that
  give a choice of extensions involved in the definition of the trivial
  Whitehead product
  $[\widetilde{E}\gamma'^1_{\lambda,1},\widetilde{E}\gamma'^1_{\lambda,1}]_w$.

  From Proposition \ref{actionmultiple}(iv) and the uniqueness up to homotopy of the wedge map that gives a Whitehead product of order $2$, it
  follows that the map $g'_{\lambda,(2)} $ has to be homotopy equivalent to
  a $N$ covering of the map $ E( \gamma_{\lambda, 1}) \vee E(  \gamma_{\lambda, 1}):S^2 \vee S^2
  \rightarrow BH_1 \subset BG^1_{\lambda} $. The latter extends  to $\C P^2 $ and give $ E_{\lambda}^{(k)}.$ 

 Therefore the restrictions of both
  $E_{\lambda}^{(k)}$ and $Q_{\lambda}^{(k)}$ to their 2-skeletons must be
  isotopic and part (i) of the lemma follows.

  Part (ii) is an immediate consequence of Theorem \ref{egwruled},
  Proposition \ref{propofpgw} and \eqref{pgwzero}.

  Part (iii) then follows from part (ii) and Corollary
  \ref{pgwobstructfibration}.
\end{proof}
\QED

 Consider now  the long exact sequences in homotopy:

\begin{math}
  \xymatrix{
    \ar[r] &
    \pi_4 \Diff_0(T^2 \times S^2) \ar@{=}[d] \ar[r] &
    \pi_4 \Aa^1_{\lambda} \ar[d] \ar[r] &
    \pi_3 G^1_{\lambda} \ar[r] &
    \pi_3 \Diff_0(T^2 \times S^2) \ar@{=}[d] \ar[r] & \\
    \ar[r] &
    \pi_4 \Diff_0(T^2 \times S^2) \ar[r] &
    \pi_4 \Aa^1_{\lambda+ t} \ar[d] \ar[r] &
    \pi_3 G^1_{\lambda +t} \ar[r] &
    \pi_3 \Diff_0(T^2 \times S^2) \ar[r] & \\
     & &  }
\end{math}\\

Take a value $\lambda $ close to $1$ and assume that $ \eta_{\lambda} $ maps
to a nontrivial element in $\pi_{3} \Diff_0(T^2 \times S^2)$. Then it follows
from Proposition \ref{mcduffdiff}(iv) that $\eta_{\alpha}$ has to lift to a
nontrivial element in $\pi_3{\cal D}_0^1$.
 But from Proposition
\ref{mcduffdiff}(ii) it follows that there should be an element $\xi$ in $G_{1}^1 $ with the same image  in  $\pi_{3} \Diff_0(T^2 \times S^2)$ as $\eta_{\lambda}$. Any extension  $ \xi_{\lambda}, \lambda>0$ must have trivial PGW in class $A-K$, from Corrolary \ref{pgwobstructfibration}. From  lemma \ref{fibersumtorus} and the symplectic sum formula it follows that we have an element $\eta'_{\lambda} $ representing the class $\widetilde{\xi}_{\lambda} -\widetilde{\eta}_{\lambda}$ with nontrivial PGW in class $A-F$ and whose image in   $\pi_{3} \Diff_0(T^2 \times S^2)$ is trivial.

Therefore $ \eta'_{\lambda} $ must be the boundary of some cycle $b_{\lambda}
:S^4 \rightarrow \Aa^1 _{\lambda}$ and this extends continuously (since $
\Aa_{\lambda} \subset \Aa_{\lambda +t}$) to a family of maps $b_{\lambda+t}
:S^4 \rightarrow \Aa_{\lambda+t}$.

The boundaries of the cycles $b_{\lambda+t}$ provide a continuous family of
maps $ \eta'_{\lambda +t}: S^3 \rightarrow G^1_{\lambda+t}, t \geq 0 $, which
have a nullhomotopic image in $ \pi_{3} \Diff_0(T^2 \times S^2)$ and hence it
does not lift to $\pi_3 {\cal D}_0^1$ for any $t \geq 0$.  Since PGW are
symplectic deformation invariants from Lemma \ref{fibersumtorus}(ii) it
follows that all elements $\eta'_{\lambda} $ must give nontrivial elements in
all $\pi_3 G_{\lambda +t}$; elements that become nullhomotopic in $\pi_3 {\cal
D}_0^1$.

But for $ t > 3/2 $ this contradicts Proposition \ref{mcduffdiff}(iv) which
says that all nontrivial elements in $\pi_3 G_{\lambda +t} $ lift to nontrivial
elements in $\pi_3 {\cal D}_0^1$. Therefore the featured Whitehead product must
be nontrivial. This concludes the proof of Proposition \ref{fragileruled}(i).

\end{proof}

\subsection{On the rational homotopy type  of {\boldmath $B\Symp_0(S^2 \times
            S^2, \omega_{\lambda})$}}\label{g=0}

The following theorem is proved in \cite{AM} and describes rational cohomology
of $G^0_{\lambda}$ and implicitly the additive structure of $\pi_*
G_{\lambda}^0 \otimes \Q$:

\begin{theorem}{(Abreu-McDuff)\cite{AM}}\label{rationaltype}
  \begin{enumerate}
    \item Let $ k <\lambda \leq k+1 $ for some natural number $k \geq 0$. We
          have
          \begin{equation}
            H^*(G_{\lambda}^0, \Q) = \Lambda (a,x,y) \otimes S(w_k)
          \end{equation}
          where $ \Lambda (a,x,y)$ is an exterior algebra with generators of
          degrees $deg \;a=1$, $deg \;x =deg\; y= 3$ and $S(w_k)$ is a
          polynomial algebra with one generator of degree $4k$.
    \item For $ k <\lambda \leq k+1 $ a complete set of generators for $\pi_*
          G_{\lambda}^0 \otimes \Q$, dual to $a,x,y,w$ respectively, is given
          by $\widetilde{\gamma}_{\lambda, 1} \in \pi_1 G_{\lambda}^0$,
          $\widetilde{\alpha}$ and $\widetilde{\eta}$ in $\pi_3 G_{\lambda}^0$
          and $\widetilde{w_k}$ in $\pi_{4k} G_{\lambda}^0$.
  \end{enumerate}
\end{theorem}

Based on the additive structure provided in Theorem \ref{rationaltype} we will
give a new proof of \ref{rationaltypeclass}.

Basically, in order to find the multiplicative structure of the ring
\eqref{classring} one must understand all the rational Samelson products among
elements in the homotopy groups $\pi_*G_{\lambda} \otimes Q$ that are
dual to the given complete set of generators $a,x,y,w$. Under suitable
conditions, each nontrivial such Samelson (hence Whitehead) product gives a
relation in the ring $H^*(BG_{\lambda}^0, \Q) $. Our contribution will be the
following:

\begin{proposition}\label{fragileruled}
  For all $k \geq 1 $ and $ k <\lambda \leq k+1 $ the Samelson product of
  order $2 k +1$, $S^{(2k+1)} (\widetilde{\gamma}'^0_{k+1, k+1}) = \{ 0,
  \widetilde{w}_k \} \subset \pi_{4k} (G^0_{k+1}) $ where $w_k$ is a
  nontrivial fragile element.
\end{proposition}

\begin{proof}{ }
As stated in Proposition \ref{mcduffacs}(ii) there exist maps
$h_{k+1,k+1+\epsilon}: G_{k+1} \rightarrow G_{k+1 + \epsilon}$ and hence one
gets maps $h'_{k+1,k+1+\epsilon}: BG_{k+1} \rightarrow BG_{k+1 + \epsilon}$.

Part (iv) in the same proposition implies that this maps induce an isomorphism
on $ \pi_i BG_{\lambda} $ for $i \leq 4k$.  In particular any map defined on a
CW-complex $B$ of dimension less or equal to $4k$, $f: B \rightarrow BG_{k+1+
\epsilon} $ must belong to a continuous family $ f_{\lambda} : B \rightarrow
BG_{\lambda}, k+1 \leq \lambda \leq k+1+\epsilon$. Let us fix $\epsilon$ small
and apply this conclusion to the map derived from (\ref{actionruledsympc}):
\begin{equation}\label{theequivmap}
  f_{k+1+\epsilon, (2k) }: \C P^{2k} \rightarrow BH_{k+1} \subset
  BG_{k+1+\epsilon}
\end{equation}
and obtain a continuous family 
\begin{equation}\label{fam}
  f'_{\lambda,(2k)} : \C P^{2k} \rightarrow B G_{\lambda}, k+1 \leq \lambda
  \leq k+1+\epsilon
\end{equation}
which for $\lambda = k+1+\epsilon$ coincides with (\ref{theequivmap}).

The maps $f'_{\lambda,(i)}: S^2 \rightarrow BG_{k+1}$, $i = 1, \ldots, 2k$
give a family of towers of fibrations of length $2k$ as in
\eqref{towerparameter}, for some elements in $\pi_2 BG_{\lambda}$ that are
homotopy equivalent to a multiple of $\widetilde{E}(\gamma'^0_{\lambda,k+1})
$. After multiplying with an appropriate high power $N$, we can assume that
the family $ \gamma'^0_{\lambda, k+1} $ satisfies the hypothesis of
Proposition \ref{Whitehead} when we set $p=2k$.

Then we claim that the fibration $ Q'^{(p)} _{k+1} \rightarrow \C P^{2k}$
cannot extend to a $ Q'^{(p+1)} _{k+1} \rightarrow \C P^{p+1}$. Indeed if it
were, then $ Q'^{(p+1)}_{k+1}$ admitted a deformation $ Q'^{(p+1)}_{\lambda}$
with respect to the parameter as in Proposition \ref{Whitehead} and from
\eqref{theequivmap} and \eqref{fam} above, it would follow that there exist an
appropriate symplectic fibration $R^{(p+1)}_{k+1+\epsilon} \rightarrow
S^{4k+2}$ such that:
\begin{equation}\label{sumfb}
  Q'^{(p+1)}_{k+1 +\epsilon} = E_{k+1+\epsilon}^{(p+1)} \#
  R^{(p+1)}_{k+1+\epsilon}
\end{equation}
for $ E_{k+1+\epsilon}^{(p+1)}$ a N-covering of the associated fibration
$Q_{k+1+\epsilon}^{(k+1),(p+1),0 }$ defined in (\ref{fibrationnot}).

As argued in the $g=1$ case, any invariants counting maps in class $A-(k+1) D$
on $Q'^{(p+1)}_{k+1 +\epsilon}$ must be zero and from Theorem \ref{egwruled}
and the symplectic sum formula (proposition \ref{propofpgw} (iii)) it follows that:

\begin{equation}\label{pgwnotzeroh}
  PGW _{1,0}^{\tilde{J} }(R^{(p+1)}_{\lambda} , s_{A-(k+1)F}) \neq 0, \lambda
  \geq k+1+\epsilon
\end{equation}

But this implies (again using Proposition \ref{mcduffacs}(iv)), that there
exist essential maps $a_{\lambda}: S^{4k+1} \rightarrow G_{\lambda}$, which is
false for sufficiently large $\lambda$.

\end{proof} 
\QED
 
\begin{remark}\label{lolo}{\rm
  We have relied throughout the paper on the classical definition of the
  higher order Whitehead products because we made use of their obstruction theoretic properties.
  Allday's \cite{A1} definition of  rational Whitehead products in the graded
  differential Lie algebra $\pi_*BG_{\lambda} \otimes \Q$ is being used in showing some of the results in the following lemma, as they sometimes use homology rather than homotopy relations. These invariants  are in one to one correspondence with the usual rational Whitehead products,
  and in this case with Samelson products in $\pi_* G_{\lambda}$.}
\end{remark}
The rational Whitehead products in $\pi_*BG_{\lambda}$ are multilinear. To
ease the computations we will write $\widetilde{A} =
\widetilde{E}(\gamma^0_{\lambda,1}) \in \pi_2 BG_{\lambda} \otimes Q$,
$\widetilde{Y} = \widetilde{E}(\alpha) \in \pi_4 BG_{\lambda} \otimes Q$,
$\widetilde{X} = \widetilde{E}(\eta) \in \pi_4 BG_{\lambda} \otimes Q$ and $
\widetilde{W_k} = \widetilde{E}(w_k) \in \pi_{4k+1} BG_{\lambda} \otimes Q$.

\begin{lemma}\label{sams}
  For any $k \geq 1 $ and  $ k <\lambda \leq k+1 $ we have 
  \begin{enumerate}
    \item Any Whitehead product of order less than $ k+1 $ is vanishing and
          also the following order $k+1$ products vanish:
          \begin{multline} \label{relationzero}
            [ \widetilde{A},\ldots , \widetilde{A},\widetilde{X} +
            \widetilde{Y}, \ldots, \widetilde{X} + \widetilde{Y}]=[
            \widetilde{A},\ldots, \widetilde{A} , \widetilde{X} + 4
            \widetilde{Y},\ldots, \widetilde{X} + 4 \widetilde{Y}] = \\
            =\ldots = [\widetilde{A},\ldots, \widetilde{A} , \widetilde{X} +
            k^2 \widetilde{Y}, \ldots, \widetilde{X} + k^2 \widetilde{Y}]=0
          \end{multline}
    \item The following Whitehead product of order $k+1$ is nontrivial and
          consists of only one element in $\pi_{4k+1} BG_{\lambda}$:
          \begin{equation}\label{thenonzerowh} 
            0 \neq [ \widetilde{A} , \widetilde{X} , \ldots, \widetilde{X}]
          \end{equation}
  \end{enumerate}
\end{lemma}
\begin{proof}{}
  Clearly $[ \widetilde{A}, \widetilde{A}]=0$. Considerations of the dimension
  of $\pi_* BG_{\lambda}$ imply that any other Whitehead products of order
  strictly less than $k+1$ must also vanish. Therefore Proposition
  \ref{arkovitz}(a) implies that any Whitehead product of order $k+1$ is
  defined and contains only one element. Since all the Lie subgroups $ K_i$,
  $i \leq k$ embed in $G_{\lambda}$, lemma \ref{actionmultiple}(i) and
  (ii) yield part (i) of the present lemma. We use here the fact that the classifying
  space of a Lie group is an H-space and hence it has vanishing rational
  Whitehead products.

  To prove the second part let us first notice that the indeterminacy in the
  Whitehead product $W^{(2k+1)} (\widetilde{A})$ obtained in Proposition
  \ref{fragileruled} implies, according with Proposition \ref{arkovitz}(a),
  that nonvanishing lower order Whitehead products must exist. Again, from
  dimension considerations, it follows that they can only be of order $p+s<
  2k+1$, $p>0, s>0$ and $ 2p+4s =4k+2$:
  \begin{equation}\label{relationabhigh}
    0 \neq [ \widetilde{A}, \ldots, \widetilde{A} , a_1\widetilde{X} + b_1
    \widetilde{Y},\ldots, a_s \widetilde{X} + b_s \widetilde{Y}]
  \end{equation} 

  \noindent\noindent{\bf Claim:} The minimum Whitehead order is $k+1$.
  
  \begin{proof}{ of the claim}
    Assume that the minimum Whitehead order is $p+s >k+1$. Hence $ p>1$ and as
    above $2p+4s =4k +2$.

    Consider the following equation in $b$:

    \begin{equation}\label{eqinb}
      0 = [ \widetilde{A}, \ldots, \widetilde{A} , \widetilde{X} + b
          \widetilde{Y},\ldots, \widetilde{X} + b \widetilde{Y}]
    \end{equation} 

    This equation has degree $s$ and coefficients in $\pi_{4k+1} BG_{\lambda}
    \otimes Q$ given by Whitehead products (containing only one element) of
    type $(p,s)$ that give a basis for all the possible Whitehead products of
    type $(p,s)$.

    Moreover, Proposition \ref{actionmultiple} implies that the equation must
    have $k$ solutions $b=1,4,\ldots k^2$ provided by the $k$ different Lie
    groups actions.

    But $k = \frac{2p+4s -2}{4} >s $ whenever $p>1$ and hence all the
    coefficients must be zero in this case. Since these coefficients generate
    all Whitehead products of the given type $(p,s)$, it follows that $p$ must
    be $1$ and hence the minimum order of an existing nontrivial product of
    type $(p,s) $ must be $k+1$.
  \end{proof}\QED

  Combined with part (ii) of our lemma this implies part (ii).
\end{proof}\QED
  
\begin{proof}{ of Theorem \ref{rationaltypeclass}}
The proof will now follow the same lines as the proof in \cite{AM}: One has to
build the Sullivan minimal model for $H^*(BG_{\lambda}, \Q)$ by giving a
complete set of generators and relations.
 
As explained in Andrews-Arkovitz a complete set of generators for the Sullivan
minimal model's differential algebra $\Mm$ of $BG_{\lambda}$ is given by
elements in the dual homotopy groups $\Hom( \pi_*(BG_{\lambda}^0 \otimes \Q,
\Q))$.

We therefore take a complete set of generators $A \in \Mm^2$, $X,Y \in \Mm^4$,
and $W_k \in \Mm^{4k+1}$, the duals of the homotopy elements
$\widetilde{A}$, $\widetilde{X}$ , $\widetilde{Y}$ and $\widetilde{W_k}$.  We
need to understand the degree 1 differential $d$ on $\Mm$.

Consider first  the case $1 < \lambda \leq 2$.

If we denote by $\Mm^{\circ}$ the quotient of $\Mm$ by the elements of degree
$0$, then any complete set of generators on $\Mm$ induces a filtration
$\Mm_s^{\circ}$ on $\Mm^{\circ}$, with ${\cal M}_s^{\circ}$ being the
subalgebra generated by products of $s$ generators.  In this case the
Whitehead minimal order is $r=2$.

According with \cite[Proposition 6.4]{AA} for any $\mu \in \Mm $ we must have
$d \mu \in \Mm_2^{\circ}$ . This, and degree considerations immediately imply
that $dA=dX=dY=0$ (all elements in $\Mm^4$ are indecomposable) and these
elements transgress to generators in $H^*(BG_{\lambda}^0, \Q)$.

Theorem 5.4 in \cite{AA} states that for any $\mu$ with $d\mu \in {\cal
M}_s^{\circ}$, and $z \in [x_1,x_2, \ldots, x_s] \in \pi_*(BG_{\lambda})
\otimes \Q$, the Sullivan pairing $\langle {\bar \mu}, z \rangle$ can be
computed in terms of suitable coefficients coming from the corresponding
universal Whitehead products. This ultimately allows one to write $d\mu$ as a
relation between the generators that will give a relation in the cohomology
ring.

All we need to find in this case is who will $d W$ be in this case. On one
hand $dW \in {\cal M}_2^{\circ}$ and on the other hand it corresponds (via
Sullivan's pairing) to a (minimal order) Whitehead product of order 2. It
follows that $dW$ must be equal to a homogeneous function $F_2$ of order two
in the remaining variables $A,X,Y$. Exactly as explained in \cite{AM}, one may
think of it as a symmetric bilinear function on a vector space spanned over
$\Q$ by the base dual to $A,X,Y$; namely, $\widetilde{A}, \widetilde{X},
\widetilde{Y}$.

But from \eqref{relationzero} and \eqref{thenonzerowh} if follows that
$F_2(\widetilde{A}, \widetilde{X}+ \widetilde{Y}) = 0$ and $F_2(\widetilde{A},
\widetilde{X}) \neq 0$ and hence $F_2 =A(X-Y)$ (up to a multiple).

The situation is similar when $k < \lambda \leq k+1$ for arbitrary $k$.

In this case the free graded differential algebra $\Mm$ has generators $A \in
\Mm^2$, $X,Y \in \Mm^4$,and $W_k \in \Mm^{4k+1}$ and the minimal Whitehead
order is $k+1$. As before, if follows that any $d \mu \in \Mm^{\circ}_{k+1} $
and hence $DA=DX=DY=0$.

In this situation $dW_k$ is in ${\cal M}_{k+1}^{\circ}$ and corresponds (via
Sullivan's pairing) to a (minimal order) Whitehead product of order $k+1$
given by (\ref{thenonzerowh}). Hence
\begin{equation}
dW=F_{k+1}(A,X,Y)
\end{equation}
where $F_{k+1} $ is homogeneous of degree $k+1$ and corresponds to a symmetric
$(k+1)$-linear function defined on a vector space spanned over $\Q$ by the
basis $\widetilde{A}, \widetilde{X}, \widetilde{Y}$.  Moreover, from Lemma
\ref{sams} we have that
\begin{equation}\label{generalf} 
  F_{k+1} (\widetilde{A},\widetilde{X} + i^2 \widetilde{Y}, \ldots,
  \widetilde{X} + i^2 \widetilde{Y})=0, \ i = 1, \ldots, k
\end{equation}
and
\begin{equation}\label{generalfnon} 
  F_{k+1} (\widetilde{A},\widetilde{X}, \ldots, \widetilde{X})\neq 0
\end{equation}

One can check that the multilinear function $F_{k+1}=A (X-Y)(X-4Y)\ldots(X-k^2
Y)$ satisfies this relations and is unique up to a constant.
\end{proof}\QED

\subsection{Higher genus cases}

\begin{proof}{ of Proposition \ref{fragileruledhg}(iii)}
  The proof will be a slightly more elaborated version of the proof of
  Proposition \ref{fragileruled}.

  Let us fix $\epsilon$ small and consider the obvious maps derived
  from \eqref{actionruledsympc}:
  \begin{equation}\label{eqmaphg}
    f_{k+\epsilon, (p)}: \C P^p \rightarrow B H_k \subset BG_{k+\epsilon}^g
  \end{equation} 

  \noindent{\bf Claim:} For $k > \lfloor g/2 \rfloor$ the map $f_{k+\epsilon,
  (g)}$ belongs to a continuous family
  \begin{equation}\label{famhg}
    f'_{\lambda, (g)}: \C P^{g} \rightarrow BG_{\lambda}^g, k \leq \lambda
    \leq k+\epsilon
  \end{equation}

  The claim follows from Proposition \ref{mcduffdiff}(ii). Indeed since the
  maps $h'_{\lambda, k+\epsilon} : BG^g_{\lambda} \rightarrow
  BG^g_{k+\epsilon}$ induce an isomorphism on $\pi_i$, $i = 1, \ldots, 2g$ for
  $\lfloor g/2 \rfloor < k \leq \lambda \leq k+\epsilon$, any map defined on a
  CW-complex of dimension less than $2g$ must extend to a continuous family as
  in \eqref{famhg}. Moreover $\pi_1 BG_{\lambda}^g \otimes Q$ is 1-dimensional
  and hence after passing to high powers (hence finite coverings of the
  induced fibration) we may assume that the family of towers of fibrations
  $Q'^{(p)}_{\lambda}$ of length $p$ given by $f^g_{\lambda,(p)}$, and their
  restrictions to the lower skeletons, gives a choice of tower for the
  Whitehead products of $\widetilde{E}(\gamma'^g_{\lambda, k}), k \leq \lambda
  \leq k+\epsilon$.

 {\bf Assumption A1}  Let us assume that all the Whitehead products $W^{(r)}
  (\widetilde{E}(\gamma'^g_{k, k})) \in \pi_{2r-1} (BG_{k}^g)$ of order $g
  \leq r \leq 2k+g-1$ vanish.

  Then the  tower $Q'^{(g)}_{k}$ of length $g$ obtained at level $\lambda = k$, must extend to
  a tower $Q'^{(2k+g-1)}_{k}$ of length $2k+g-1$ as in Proposition
  \ref{Whitehead}point(b).
 Furthermore, the tower extends again with respect to the
  parameter and we obtain once more families of fibrations
  \begin{equation}
    Q'^{(r)}_{\lambda} {\rm \ for \ } 1 \leq r \leq 2k+g-1 {\rm \ and \ } k
    \leq \lambda \leq k+\epsilon
  \end{equation} 

As before, Denote by $E_{k+\epsilon}^{(2k+g-1)}$  a $N$-covering of the fibration \eqref{fibrationnot} arising from the circle action $H_k$.

The construction above implies that two fibration  $Q'^{(2k+g-1)}_{k+\epsilon} \rightarrow \C P^{2k+g-1}$ and $ E_{k+\epsilon}^{(2k+g-1)} \rightarrow \C P^{2k+g-1}$ agree over the $2g$ skeleton.
 We should point out that the corresponding restriction to the $2g$ skeletons are just the fibrations  $Q'^{(g)}_{k+\epsilon} \rightarrow \C P^{g}$ and $ E_{k+\epsilon}^{(g)}$ in the towers. 
Under the vanishing assumption A1 we have the following:
\begin{lemma}\label{fis}
There exist fiberwise symplectic deformations  $Q'^{(2k+g-1)}_{\lambda} \rightarrow \C P^{2k+g-1}$ and \\$ E_{\lambda}^{(2k+g-1)} \rightarrow \C P^{2k+g-1}$ such that, for a value $\lambda=a$ sufficiently large, the two corresponding fibrations are symplectically isotopic.
\end{lemma}
\begin{proof}{}
Let us remind the reader that $\pi_i B{\cal D}_0^g $  is trivial if $i>4$. That will definitely be the case when $i>2g\geq4$. 

Firstly, let us notice that since they agree (isotopic would be enough)  on the the $2g$ skeleton we have
  \begin{equation}
    Q'^{(g+1)}_{k+\epsilon}= E_{k+\epsilon}^{(g+1)} \# R^{(g+1)} _{k+\epsilon} \end{equation}
where  $R^{g+1} _{k+\epsilon} \rightarrow S^{(2g+1)} $ is a symplectic
  fibration that corresponds to an element $\tilde{\eta}_{k+\epsilon}$ in  $\pi_{2g+1} BG^g_{k+\epsilon}$. Since  (since $2g+1>4$) the image of  $\tilde{\eta}_{k+\epsilon}$ in   $\pi_i B{\cal D}_0^g $ vanishes and hence so does its image in $\pi_{2g+1}  B{\cal D}_0^g $. THerefore $\tilde{\eta}_{k+\epsilon}$ hence it must be the image of an element $\tilde{\beta}_{k+\epsilon}\in \pi_{2g+2} \Aa^g_{k+ \epsilon}$  in the following diagram:

 Consider now  the long exact sequences in homotopy:

\begin{math}
  \xymatrix{
    \ar[r] &
    \pi_{2g+2} \Diff_0(\Sigma_g^2 \times S^2) \ar@{=}[d] \ar[r] &
    \pi_{2g+2} \Aa^g_{k+ \epsilon} \ar[d] \ar[r] &
    \pi_{2g+1} G^g_{k+ \epsilon} \ar[r] &
    \pi_{2g+1} \Diff_0(\Sigma_g \times S^2) \ar@{=}[d] \ar[r] & \\
    \ar[r] &
    \pi_{2g+2} \Diff_0(\Sigma_g \times S^2) \ar[r] &
    \pi_{2g+2} \Aa^g_{\lambda} \ar[d] \ar[r] &
    \pi_{2g+1} G^g_{\lambda} \ar[r] &
    \pi_ {2g+1}\Diff_0(\Sigma_g^2 \times S^2) \ar[r] & \\
     & &  }
\end{math}\\

 There clearly is a continuous family  $\beta_{\lambda} :S^{2g+2} \rightarrow \Aa^g_{\lambda}$ whose boundaries give a continuous extension $\eta_{\lambda}: S^{2g+1} \rightarrow BG_{\lambda}^g$, $\lambda \geq k+ \epsilon$.
 Assume that the latter family  is essential in the homotopy of $BG^g_{\lambda}$ for all values of $\lambda$. Then  $\tilde{\beta}_{\lambda}\in \pi_{2g+2} \Aa^g_{\lambda}$ must be nontrivial for all $\lambda$, and it immediately follows that it must also give essential element $\beta_{\infty}$ in
$\pi_{2g+2} \Aa^g_{\infty}$, where $\Aa^g_{\infty} =\cup_{\lambda >0} \Aa^g_{\lambda}$.
We now look at the homotopy fibration:
$$G^g_{\infty} \rightarrow{\rm Diff}_0( \Sigma_g \times S^2) \rightarrow \Aa^g_{\infty} $$ 
In the long exact sequence in homotopy for this fibration the element   $\beta_{\infty}$ does not lift to $\pi_{2g+2}\Diff_0(\Sigma_g \times S^2)$. This is because none of the elements $\tilde{\beta}^{2g+1}_{\lambda}$  lift to $\pi_{2g+2}\Diff_0(\Sigma_g \times S^2)$. Therefore $\beta_{\infty}$ must have a nontrivial boundary and give an nontrivial element $\eta_{\infty } \in \pi_{2g+1} G^g_{\infty}$ which is impossible from Proposition \ref{mcduffdiff} part(ii). 

Therefore  there exist   a value $\lambda_0$ where $\beta^{2g+1}_{\lambda_0} $ becomes inessential in $BG_{\lambda_0}$.
Using the continuous family  $\beta_{\lambda} $ we obtain two fiberwise deformations of   $ Q'^{(2g+1)}_{\lambda}$ and $ E_{\lambda}^{(2g+1)} $, $\lambda \geq k+ \epsilon$, that are isotopic at $\lambda = \lambda_0$. 
On them we build two deformations  $Q'^{(2k +g -1)}_{\lambda}$ and
$ E_{\lambda}^{(2k +g -1)}$, $\lambda \ge k+ \epsilon$, whose restrictions to the $2g+1$ skeleton are isotopic at $\lambda=\lambda_0$. 

We repeat this process $ 2k$ more steps, in each of which we obtain fibrations that agree on a skeleton of dimension two bigger.
In the end we get a large  value $\lambda=a$ and two fiberwise deformations  $Q'^{(2k+g-1)}_{\lambda}$ and $ E_{\lambda}^{(2k+g-1)}$, $k+ \epsilon \leq \lambda \leq a$ such that the following isotopy holds 
\begin{equation}
 Q'^{(2k+g-1)}_{a}  \approx  E_{a}^{(2k+g-1)} \end{equation}
\end{proof}

From here the result is straightforward: on one hand  $Q'^{(2k+g-1)}_{\lambda}$ must have trivial PGW in class $A-kF$ due to corollary \ref{pgwobstructfibration} and on the other hand $E_{\lambda}^{(2k+g-1)}$ must have a nontrivial PGW in class  $A-kF$ from Theorem \ref{egwruled} and invariance under deformation.

But this is impossible, therefore the assumption A1 must be false.
\end{proof}

\begin{remark} 
\begin{enumerate}
  \item It is very likely that the proposition \ref{fragileruledhg} (iii)  can be strengthened to a statement about
  finding the nontrivial elements in the Samelson products $ S^{(2k+g-1)}
  (\widetilde{\gamma'}^g_{k, k}) \in \pi_{4k+2g}G_{k}^g$, as conjectured in
  \cite{M}.  For that, we need to extend the gluing argument from
  \cite{M} page 20 to show that the maps         
  \begin{equation*} f_{k+\epsilon, (2k+g-2)} : \C P^{2k+g-2} \rightarrow BH_{k} \subset
    BG^g_{k+\epsilon}
  \end{equation*}
  extend to continuous families as in \eqref{famhg}.
\item The lemma \ref{fis} proves in fact that any two fibrations that agree on a the $4$ skeletons deform into isotopic fibrations for large $\lambda$. 
\end{enumerate}

\end{remark}

\noindent{\bf Acknowledgments} This paper owes enormously to my former advisor
Dusa McDuff. I would like to thank her for suggesting the original problem, as
well as for her interest and suggestions. I would like to thank Tom Parker for
useful discussions.


\begin{thebibliography}{99999}

\bibitem{A} M. Abreu,  Topology of symplectomorphism groups of $S^2\times
  S^2$, {\sf Inv. Math.}, {\bf 131} (1998), 1-23.

\bibitem{AM} M. Abreu and D. McDuff, Topology of symplectomorphism
groups of rational ruled surfaces, {\it J. Amer. Math. Soc.}  {\bf
13}, no.4 (2000), 971-1009.


\bibitem{AA} P. Andrews and M. Arkovitz, Sullivan's minimal model and higher
order Whitehead products, {\it Can. J. Math} {\bf 13} (1978), 991-982.

\bibitem{A1} C. Allday, Rational Whitehead products and a spectral sequence of
Quillen,II, {\it Hust. J. Math.}  {\bf 3} no. 3 (1977), 301-308.

\bibitem{A2} C. Allday, Rational Whitehead products and a spectral sequence of
Quillen, {\it Pac. J. Math.}  {\bf 46} no. 2 (1973), 313-323.

\bibitem{Au} M. Audin, Hamiltonians on four-dimensional compact symplectic manifolds, Geometrie Symplectique et mecanique, Lect. Note in Amth., 1416, Springer, Berlin (1996).

\bibitem{B} O. Bu\c se, Parametric Gromov-Witten invariants and
symplectomorphism groups, to appear in {\it Pac. J. Math.}.

\bibitem{G} M. Gromov, Pseudo holomorphic curves in symplectic manifolds,
            {\it Inv. Math.}, {\bf 82} (1985), 307-347.

\bibitem{GS} V. Guillemin, E. Lerman, S. Sternberg, Symplectic fibrations and multiplicity diagrams, Combridge Univ. Press,Combridge, 1996.

\bibitem{K}  Y. Karshon, Periodic hamiltonian flows on four-dimensional manifolds,Mem.Amer.Math.Soc.141, no. 672,(1999) (viii+71 pp).

\bibitem{LO} H.V. Le and K. Ono Parameterized Gromov-Witten and topology of
symplectomorphism groups, preprint, (2001).

\bibitem{LT} J. Li and G. Tian Virtual moduli cycles and Gromov-Witten invariants of general symplectic manifolds, Topics in $4$-manifolds,(Irvine,CA,1996),47-83,First Int Press Lect. Ser.I, Internat. press (1996).

\bibitem{M} D. McDuff,  Symplectomorphism groups and almost complex
  structures, Essays on geometry and related topics, Vol 1,2, 527-556, {\sf Enseignement Math.,Geneva}(2001).

\bibitem{MS1} D. McDuff and D.A. Salamon,  $J$-holomorphic curves and
quantum cohomology, University Lecture Series {\bf 6}, American Mathematical
Society, Providence, RI, 1994.

\bibitem{MS2} D. McDuff and D.A. Salamon,  Introduction to Symplectic
Topology, 2nd edition, Oxford University Press, 1998.

\bibitem{MT} D. McDuff and S.Tolman, Topological properties of
hamiltonian circle actions, preprint( 2004)SG/0404340 .

\bibitem{P} G. Porter, Higher order Whitehead products, {\it Topology}  {\bf
3}  (1965), 123-135.

\bibitem{S1} P. Seidel, Floer homology and the symplectic isotopy problem,
Ph.D.  thesis, Oxford University, 1997.

\bibitem{S2} P. Seidel, $\pi_1$ of symplectic automorphism groups and
invertibles in quantum homology rings, Geon. Func. Anal. 7(1997), no. 6,
1046-1095, Oxford University, 1997.

\bibitem{Sh} G.-Y. Shi, in preparation.

\end{thebibliography}
\end{document}